\documentclass[12pt]{amsart}

\usepackage[text={6.5in,9in},centering]{geometry}

\usepackage{amsmath,amssymb,amsfonts}
\usepackage{graphicx}
\usepackage{bm}
\usepackage{color}
\usepackage{fix-cm}
\usepackage{latexsym,mathrsfs}

\usepackage{color}
\usepackage[caption=false,labelfont={rm,rm},textfont=rm]{subfig}

\usepackage[english]{babel}
\graphicspath{{FigPDF/}{FigEPS/}}

\newtheorem{theorem}{Theorem}[section]
\newtheorem{lemma}[theorem]{Lemma}
\newtheorem{proposition}[theorem]{Proposition}

\newtheorem{algorithm}[theorem]{Algorithm}

\newtheorem{remark}[theorem]{Remark}

\numberwithin{equation}{section}
\numberwithin{figure}{section}
\numberwithin{table}{section}
\renewcommand{\arraystretch}{1.5}

\newcommand{\vect}[1]{\bm{#1}}

\newcommand{\K}{T} 
 
\renewcommand{\O}{\Omega} 
 
\renewcommand{\div}{{\rm div}} 
\newcommand{\grad}{{\boldsymbol \nabla}}

\newcommand{\n}{{\bf n}} 
\newcommand{\uh}{\vect{u}_h}

\newcommand{\jump}[1]{\lbrack\!\lbrack #1 \rbrack\!\rbrack}

\newcommand{\Eh}{{\mathcal{E}_h}} 
\newcommand{\PzeroE}{\mathcal{P}^0_{E}}

\newcommand{\Eho}{{\mathcal{E}^{o}_h}} 
\newcommand{\EhN}{{\mathcal{E}^{N}_h}} 
\newcommand{\EhD}{{\mathcal{E}^{D}_h}} 
\newcommand{\Ehb}{{\mathcal{E}^{\partial}_h}}

\renewcommand{\lor }{\longrightarrow}

\newcommand{\Th}{\mathcal{T}_h} 
\newcommand{\hK}{\h_{\K}} 
\newcommand{\h}{h}

\renewcommand{\div}{\textrm{div}}

\newcommand{\calB}{\mathcal{B}}
\newcommand{\calA}{ \mathcal{A}}

\newcommand{\calP}{ \mathcal{P}}

\newcommand{\calE}{ \mathcal{E}}
\newcommand{\calZ}{ \mathcal{Z}}

\newcommand{\bu}{\vect{u}}
\numberwithin{equation}{section}
\renewcommand{\arraystretch}{1.5}

\newcommand{\Reals}[1]{{\rm I\! R}^{#1}}
\newcommand{\Space}[2]{#1^{\textrm{#2}}}
\newcommand{\Spacev}[2]{\vect{#1}^{\textrm{#2}}}
\newcommand{\jumpv}[1]{\lbrack\!\lbrack\vect{#1}\rbrack\!\rbrack} 
\newcommand{\avg}[1]{\{\!\!\{#1\}\!\!\}} 
\newcommand{\avgv}[1]{\{\!\!\{\vect{#1}\}\!\!\}} 

\newcommand{\scalarT}[2]{(#1,#2)_{\mathcal{T}_h}}
\newcommand{\scalarTT}[2]{(#1:#2)_{\mathcal{T}_h}}

\newcommand{\de}{{:=}}

\newcommand{\Gsym}[1]{\vect{\varepsilon}(\vect{#1})}
\newcommand{\Ceps}[1]{\mathcal{C}\vect{\varepsilon}(\vect{#1})}
\newcommand{\eps}[1]{\vect{\varepsilon}(\vect{#1})}

\newcommand{\triplenorm}[1]{
  \left\vert\kern-0.9pt\left\vert\kern-0.9pt\left\vert #1
  \right\vert\kern-0.9pt\right\vert\kern-0.9pt\right\vert}  
\title[A subspace correction for dG discretizations of elasticity]{A subspace correction method for ciscontinuous
 Galerkin discretizations of linear elasticity equations}
\begin{document}

\author{Blanca Ayuso de Dios}\address{Centre de Recerca Matem\`{a}tica,
  Campus de Bellaterra, 08193 Bellaterra, (Barcelona), Spain. Email:
  \texttt{bayuso@crm.cat}}

\author{Ivan Georgiev}\address{
  Johann Radon Institute for
  Computational and Applied Mathematics, Austrian Academy of
  Sciences Altenberger Str. 69, 4040 Linz,
  Austria. Email: \texttt{ivan.georgiev@oeaw.ac.at}
}
\author{Johannes Kraus}\address{Johann Radon Institute for
  Computational and Applied Mathematics, Austrian Academy of
  Sciences Altenberger Str. 69, 4040 Linz,
  Austria. Email:\texttt{johannes.kraus@oeaw.ac.at}
}
\author{Ludmil Zikatanov}\address{Department of Mathematics, The
  Pennsylvania State University, University Park, PA 16802, USA.
Email: \texttt{ltz@math.psu.edu}
}
\date{\today}
\begin{abstract}
  We study preconditioning techniques for discontinuous Galerkin
  discretizations of isotropic linear elasticity problems in primal
  (displacement) formulation. We propose subspace correction methods
  based on a splitting of the vector valued piecewise linear
  discontinuous finite element space, that are optimal with respect to
  the mesh size and the Lam\'e parameters. The pure displacement, the
  mixed and the traction free problems are discussed in detail.  We
  present a convergence analysis of the proposed preconditioners and
  include numerical examples that validate the theory and assess the
  performance of the preconditioners.
\end{abstract}



\maketitle

\section{Introduction}
The finite element approximation of the equations of isotropic linear
elasticity may be accomplished in various ways. The most
straightforward approach is to use the primal formulation and
conforming finite elements. It is well known that such a method, in
general, does not provide approximation to the displacement field when
the material is nearly incompressible (the Poisson ratio is close to
$1/2$). This phenomenon is called \emph{volume locking}. To alleviate
locking, several approaches exist. Among the possible solutions, we
mention the use of mixed methods, reduced integration techniques,
stabilization techniques, nonconforming methods, and the use of
discontinuous Galerkin methods.  We refer to \cite{FalkR-1991aa,
  HansboP_LarsonM-2002aa} for further discussions on such difficulties
and their remedies.  In this work we focus on the Symmetric Interior
Penalty discontinuous Galerkin (SIPG) methods introduced in
\cite{HansboP_LarsonM-2002aa, HansboP_LarsonM-2003aa, WihlerT-2004aa,
  WihlerT-2006aa} for the approximation of isotropic linear
elasticity. 
We have chosen to work with these DG discretizations, since we have in
mind a method that is simple but still applicable to different types
of boundary conditions. In fact, unlike classical low order
non-conforming methods (see \cite{FalkR-1991aa}), the Interior Penalty
(IP) stabilization methods introduced in \cite{HansboP_LarsonM-2002aa,
  HansboP_LarsonM-2003aa} can be shown to be stable in the case of
essential (Dirichlet or pure displacement) boundary conditions, or
natural (Neumann type, or traction free) boundary conditions. As a
consequence, these IP methods provide a robust approximation to the
displacement field and avoid the volume locking regardless the
boundary conditions of the problem.

For the design of the preconditioners we follow the ideas introduced
in \cite{Ayuso-de-DiosB_ZikatanovL-2009aa} for second order elliptic
problems. 
However, such extensions are not straightforward, since we aim at
constructing preconditioners that work well for three different types
of boundary conditions: essential, natural and mixed boundary
conditions, used in linear elasticity. This complicates the matters
quite a bit.  We consider a splitting of the vector valued, piecewise
linear, discontinuous finite element space, into two subspaces: the
vector valued Crouzeix-Raviart space and a space complementary to it
which consists of functions whose averages are $L^2$ orthogonal to the
constants on every edge/face of the partition. This space
decomposition is direct and the spaces are orthogonal with respect to
a bilinear form obtained via using ``reduced integration'' to
calculate the contributions of the penalty terms in SIPG.

In the pure displacement case (essential boundary conditions), the
restriction of the bilinear form based on reduced integration is
coercive on the Crouzeix-Raviart space and is spectrally equivalent to
the SIPG bilinear form.  The space decomposition mentioned above is
then orthogonal in this reduced integration bilinear form.  Thus, in
case of essential boundary conditions we have a natural block diagonal
preconditioner for the linear elasticity problem: (1) a solution of a
problem arising from discretization by nonconforming Crouzeix-Raviart
elements; (2) solution of a well-conditioned problem on the
complementary space.

For traction free problems or problems with Dirichlet conditions only
on part of the boundary, the situation is quite different. On one hand
the reduced integration bilinear form when restricted to the
Crouzeix-Raviart space has a null space whose dimension depends on the
size of the problem (see \cite{FalkR-1991aa}). On the other hand in
the full SIPG bilinear form (without reduced integration) the space
splitting discussed above is no longer orthogonal. Our approach in
resolving these issues is based on a delicate estimate given in
\S\ref{subsect:cbs} which shows a uniform bound on the angle between
the Crouzeix-Raviart and its complementary space in the SIPG
bilinear form for all types of boundary conditions. Once such a bound
is available we show that a uniform block diagonal preconditioner can
be constructed.

The rest of the paper is organized as follows.  We present the linear
elasticity problem, the basic notation and discuss the DG
discretizations considered in~\S\ref{DG-IP-discretizations}.  Next,
in~\S\ref{sect:split} we introduce the splitting of the vector valued
piecewise linear DG space and discuss some properties of the related
subspaces. In section~\S\ref{sect:preconditioning}, we introduce the
subspace correction methods, and we prove that they give rise to a
uniform preconditioner for the symmetric IP method. The last
section~\S\ref{sect:numerics} contains several numerical tests that
support the theoretical results.

\section{Interior Penalty Discontinuous Galerkin methods for linear  elasticity equations}\label{DG-IP-discretizations}

In this section, we introduce the linear elasticity problem together
with the basic notation and the derivation of the Interior Penalty
(IP) methods and we discuss the stability of these methods.

\subsection{Linear Elasticity: Problem formulation and notation}
\label{sect:problem-formulation-and-notation}
Let $\Omega\subset \Reals{d}$, $d=2,3$, be a polygon or polyhedron
(not necessarily convex) and let $\vect{u}$ be a vector field in
$\Reals{d}$, defined on $\Omega$ such that $\vect{u}\in
[H^1(\Omega)]^d$.  The elasticity tensor, which we denote by
$\mathcal{C}$, is a linear operator, i.e., $\mathcal{C}:
\Reals{d\times d}_{\textrm{sym}} \mapsto \Reals{d\times
  d}_{\textrm{sym}}$, acting on a symmetric matrix $A\in
\Reals{d\times d}_{\textrm{sym}}$, in the following way:
\[
\mathcal{C}\;A = 2\mu A + \lambda \operatorname{trace}(A) I,
\]
where $\mu $ and $\lambda$ are the Lam\'e parameters and satisfy
$0<\mu_1<\mu<\mu_2$ and $0\leq \lambda <\infty$. In terms of the
modulus of elasticity (Young's modulus), $\mathfrak{E}$, and Poisson's
ratio, $\nu$, the Lam\`e parameters can be rewritten in the case of
plane strain as: $\mu=\mathfrak{E}/(2(1+\nu)) $ and
$\lambda=\nu\mathfrak{E}/ ((1+\nu)(1-2\nu)$. The material tends to the
incompressible limit (becomes incompressible) when the Lam\'e
parameter $\lambda \to \infty$ or equivalently when the Poisson's
ratio $\nu\to 1/2$.

One can show that the linear operator $\mathcal{C}$ is selfadjoint and
has two eigenvalues: (1) a simple eigenvalue equal to $(2\mu +
d\lambda)$ corresponding to the identity matrix; (2) an eigenvalue
equal to $2\mu$, corresponding to the $\frac{d(d+1)}{2}-1$ dimensional space
of traceless, symmetric, real matrices.  Thus for $d=2,3$, we always have that
\begin{equation}\label{eigsC}
2\mu\langle A:A\rangle 
\le \langle \mathcal{C}A:A\rangle \le (2\mu + d\lambda)\langle
A:A\rangle ,
\end{equation}
where $\langle\cdot : \cdot\rangle$  denotes  the Frobenius scalar product of two tensors in $\Reals{d\times d}$. We also denote by $\langle\cdot,\cdot\rangle$  the Euclidean scalar product of two vectors  in $\Reals{d}$, i.e.,
\[
\langle \vect{v},\vect{w}\rangle = \sum_{k=1}^d v_k w_k, \qquad
\langle \vect{v}:\vect{w}\rangle = \sum_{j=1}^d\sum_{k=1}^d
v_{jk} w_{jk}.
\]
The corresponding inner products in $[L^2(\O)]^d$ and  $[L^2(\O)]^{d\times
  d}$ are denoted by
\[
(\vect{v},\vect{w}) = \int_{\O}\langle\vect{v},\vect{w}\rangle,
\qquad
(\vect{v}:\vect{w}) = \int_{\O}\langle\vect{v}:\vect{w}\rangle .
\]

We write $\partial\Omega=\Gamma_N\cup \Gamma_D$ with $\Gamma_N$ and
$\Gamma_D$ referring respectively to the subsets of the $\partial \O$
where Neumann and Dirichlet boundary conditions are imposed.

Let $\Gsym{u}=\frac{1}{2} ( \grad \bu + (\grad \bu)^{T})$ be the
symmetric part of the gradient of a vector valued function $\vect{u}$. The
elasticity problem in primal formulation then is: 
Find $\vect{u} \in [H^{1+\alpha}_{\Gamma_D}(\Omega)]^d$, $\alpha>0$, 
which is the unique minimizer of the energy
functional $\mathcal{J}(\vect{u})$, given by
\begin{equation}\label{eq:minimize0}
\mathcal{J}(\vect{u}):=
\frac12
(\Ceps{u}:\Gsym{u})- (\vect{f},\vect{u}) -(\vect{g}_N, \vect{v})_{\Gamma_N} 
\end{equation}
Here $\vect{f}\in [L^2(\Omega)]^d$ is a given volume force and
$\vect{g}_N \in [H^{3/2}(\Gamma_N)]^d$ is a given surface force acting
on $\Gamma_N\subset \partial\O$. The Euler-Lagrange equations
corresponding to the minimization problem~\eqref{eq:minimize0}
give the following well known system of linear PDEs for the
unknown displacement field $\vect{u}$:
\begin{equation}\label{eq:PDEform}
\begin{aligned}
-\div (\mathcal{C}\Gsym{u}) & =\vect{f}, \quad &&\mbox{  on  } \O, \\ 
(\mathcal{C}\Gsym{u})\vect{n} & =\vect{g}_{N}, \quad &&\mbox{  on  }\Gamma_N,\\
\vect{u}&=\vect{0}, \quad && \mbox{  on  } \Gamma_D. 
\end{aligned}
\end{equation}
In the above equations, $\vect{n}$ is the outward unit normal vector
to $\partial\O$. The solution $\vect{u}$ vanishes on a closed part of
the boundary $\Gamma_D$ (Dirichlet boundary) and the normal stresses
are prescribed on $\Gamma_N$ (Neumann part of the boundary).  In the
traction free case ($\Gamma_N=\partial\Omega$), the existence of a
unique solution to \eqref{eq:PDEform} is guaranteed if the data
satisfy the following compatibility condition:
\begin{equation*}
\int_{\Omega} \vect{f}\cdot \vect{v} dx+\int_{\partial \Omega} \vect{g}_{N}\cdot \vect{v} ds =0 \quad \forall\, \vect{v} \in \mathbf{RM}(\O),
\end{equation*}
where $ \mathbf{RM}(\O)$ is the space of rigid motions, defined by:
\begin{equation}\label{rm00}
 \mathbf{RM}(\O):=\left\{ \vect{v}=\vect{a}+\vect{b}\vect{x} \quad :\quad \vect{a}\in \mathbb{R}^{d} \quad \vect{b}\in \mathfrak{so}(d)\,\,\right\}
 \end{equation}
 where $\vect{x}$ is the position vector function in $\Omega$ and $
 \mathfrak{so}(d)$ is the Lie algebra of skew-symmetric $d\times d$
 matrices. In this case, the uniqueness of solution is guaranteed up
 to a rigid motion (and is unique, if we require that the solution is
 orthogonal to any element from $\mathbf{RM}(\Omega)$).  In the case
 of $\Gamma_D\neq \emptyset$ and closed with respect to
 $\partial\Omega$ no extra conditions are required to guarantee
 uniqueness. By considering the variational formulation of
 \eqref{eq:PDEform}, the issue of solvability and uniqueness of the
 problem reduces to show coercivity of the associated bilinear
 form. As it is well known, for linear elasticity, this hinges on the
 classical Korn's inequality \cite{DuvautG_LionsJL-1976aa} which
 guarantees the existence of a generic positive constant
 $C_{\Omega}>0$ such that:
\begin{equation}\label{korn:cont}
\|\nabla \vect{v}\|_{0,\Omega}^{2} \leq C_{\Omega}\left( \|\Gsym{v}\|_{0,\Omega}^{2} +\|\vect{v}\|_{0,\Omega}^{2}\right), \qquad \forall\, \vect{v}\in [H^{1}(\Omega)]^{d}\;.
\end{equation}
The second term on the right hand side can be omitted as follows from
the Poincar\'e or Poincar\'e-Friedrich's inequality, obtaining thus
first Korn's inequality for $\vect{v}\in
[H^{1}_{0,\Gamma_D}(\Omega)]^{d}$ and second Korn's inequality for
$\vect{v}\in [H^{1}(\Omega)]^{d}/\mathbf{RM}(\Omega)$.

\subsection{Interior penalty methods: Preliminaries and notation}
We now introduce the basic notations and tools needed for the
derivation of the DG methods.

\noindent {\bf Domain partitioning.}  Let $\Th$ be a shape-regular
of partition of $\O$ into $d$-dimensional simplices $\K$ (triangles
if $d=2$ and tetrahedrons if $d=3$). We denote by $\hK$ the diameter
of $\K$ and we set $\h=\max_{\K \in \Th} \hK$. We also assume that
$\Th$ is conforming in the sense that it does not contain hanging nodes.
A face (shared by two neighboring elements or being part of the
boundary) is denoted by $E$. Clearly, such a face is a
$(d-1)$ dimensional simplex, that is, a line segment in two dimensions
and a triangle in three dimensions.  We denote the set of all faces
by $\Eh$, and the collection of
all interior faces and boundary faces by $\Eho$ and $\Ehb$,
respectively. Further, the set of Dirichlet faces is denoted by
$\EhD$, and the set of Neumann faces  by $\EhN$.
We thus have, 
\[
\Eh=\Eho\cup\Ehb,\quad 
\EhD=\Ehb\cap \Gamma_D,\quad 
\EhN=\Ehb\cap \Gamma_N,\quad
\Ehb=\EhD\cup \EhN.
\]

\noindent {\bf Trace operators (average and jump) on $E\in\Eh$.} \label{subsect:traces}
To define the average and jump trace operators for an interior face
$E\in \Eho$, and any $\K\in \Th$, such that $E\in \partial \K$ we set
$\vect{n}_{E,T}$ to be the unit outward (with respect to $T$) normal
vector to $E$.  With every face $E\in \Eho$ we also associate a unit
vector $\vect{n}_E$ which is orthogonal to the $(d-1)$ dimensional affine
variety (line in 2D and plane in 3D) containing the face.  For the
boundary faces, we always set $\vect{n}_E=\vect{n}_{E,T}$, where $T$
is the \emph{unique} element for which we have $E\subset \partial T$.
In our setting, for the interior faces, the particular direction of
$\vect{n}_E$ is not important, although it is important that this
direction is fixed.  For every face $E\in \Eh$, we define $T^{+}(E)$
and $T^{-}(E)$ as follows:
\begin{equation}\label{eq:t-plus-minus}
\begin{array}{rcl}
T^{+}(E)&\de& \{T\in \Th\  \mbox{such that} \ E\subset \partial T, \ 
\mbox{and} \ \langle\vect{n}_{E}, \vect{n}_{E,T}\rangle>0\}, \\
T^{-}(E)&\de& \{T\in \Th\ \mbox{such that}  \ E\subset \partial T, \
\mbox{and} \ \langle\vect{n}_{E}, \vect{n}_{E,T}\rangle < 0\}.
\end{array}
\end{equation}
It is immediate to see that both sets defined above contain \emph{no
  more than} one element, that is: for every face we have exactly one
$T^{+}(E)$ and for the interior faces we also have exactly one
$T^{-}(E)$. For the boundary faces we only have $T^{+}(E)$. In the
following, we write $T^\pm$ instead of $T^\pm(E)$, when
this does not cause confusion and ambiguity.

For a given function $\vect{w}\in [L^2(\O)]^d$ the average and jump trace operators
for a fixed $E\in \Eho$ are as follows:
\begin{equation}\label{eq:definition0}
\avg{\vect{w}}:=\left(\frac{\vect{w}^+ + \vect{w}^-}{2}\right),\qquad
\jump{\vect{w}}:= (\vect{w}^+ - \vect{w}^-),
\end{equation}
where $\vect{w}^{+}$ and $\vect{w}^{-}$ denote respectively, the
traces of $\vect{w}$ onto $E$ taken from within the interior of
$\K^{+}$ and $\K^{-}$.  On boundary faces $E\in \Ehb$, we set
$\avg{\vect{w}}=\vect{w}$ and $\jump{\vect{w}}=\vect{w}$.  We remark
that our notation differs from the one used in
\cite{ArnoldD_BrezziF_CockburnB_MariniL-2001aa},
\cite{ArnoldD_BrezziF_MariniL-2005aa},
\cite{ArnoldD_BrezziF_FalkR_MariniL-2007aa} (which is considered a
classical one for the IP methods).  We have chosen a notation that is
consistent with the one used in \cite{HansboP_LarsonM-2003aa}, where
the IP method we consider was introduced for the pure
displacement problem. In addition, it seems that such a choice leads to
a shorter and simpler description of the preconditioners we propose here.

\noindent {\bf Finite Element Spaces.} The piecewise linear DG space  is defined by
\begin{equation*}
\Space{V}{DG}:=\{u\in L^{2}(\O)~\mbox{such that}~u\big|_{\K} \in \mathbb{P}^{1}(\K), \quad \forall\, \K\in \Th\,\},
\end{equation*}
where $\mathbb{P}^{1}(\K)$ is the space of linear polynomials on
$\K$. The corresponding space of vector valued functions is defined as
\begin{equation*}
\Spacev{V}{DG}:=[\Space{V}{DG}]^d.
\end{equation*}
For a given face $E$, we denote by $\PzeroE :L^{2}(E)\mapsto
\mathbb{P}^{0}(E)$ the $L^{2}$-projection onto the constant
(vector valued or scalar valued) functions on $E$ defined by
\begin{align}
\PzeroE w = \frac{1}{|E|}\int_E w &\quad\mbox{for all}\quad
w\in L^{2}(E), \label{eq:l2edge}\\
\PzeroE \bm{w} = \frac{1}{|E|}\int_E \bm{w} &\quad\mbox{for all}\quad
\bm{w}\in [L^{2}(E)]^d. \label{eq:l2edgev}
\end{align}
Observe that for $\bm{w}\in \Spacev{V}{DG}$ the mid-point integration
rule implies that $\PzeroE \bm{w}=\bm{w}(m_E)$ for all $E\in
\Eh$, with $m_E$ denoting the barycenter of the edge or face $E$.

The classical Crouzeix-Raviart finite element
space can be defined as a subspace of $V^{DG}$, as follows:
\begin{equation}\label{defCRs}
\Space{V}{CR}=\left\{ v\in \Space{V}{DG} \,  \,  \,: \quad \PzeroE 
\jump{v}=0, \,\, \forall\, E\in \Eho\right\}.
\end{equation}
The corresponding space of vector valued functions is
\begin{equation}\label{defCR}
\Spacev{V}{CR}:=[\Space{V}{CR}]^d
\end{equation}

\subsection{Weighted residual derivation of the IP methods}
In \cite{HansboP_LarsonM-2003aa} the authors introduced a symmetric
interior penalty method for the problem of linear elasticity
\eqref{eq:PDEform} in the pure displacement case (i.e,
$\Gamma_D=\partial\O, \,\,\Gamma_N=\emptyset$). We define the function
space
\begin{equation*} 
[H^{2}(\Th)]^{d}=\left\{ \vect{u} \in
    [L^{2}(\O)]^{d}~\mbox{such that}~\vect{u}\big|_{\K} \in
    [H^{2}(\K)]^{d}, \quad \forall\, \K\in \Th\,\right\}.
\end{equation*}
For any pair of vector fields (or tensors)
 $\vect{v}$ and $\vect{w}$, we denote
\[
\scalarT{\vect{v}}{\vect{w}} = \sum_{\K\in \Th}\int_{\K}\langle \vect{v},\vect{w}\rangle.
\]
For scalar and vector valued functions we also use the notation
\begin{equation}\label{eq:def-edge-product}
(v,w)_{\mathcal{E}} = \sum_{E\in \calE}\int_E vw,
\quad\mbox{and}\quad
({\vect{v}},{\vect{w}})_{\mathcal{E}} = \sum_{E\in \calE}\int_E\langle
\vect{v},\vect{w}\rangle\;.
\end{equation}

We now derive, using the weighted residual
framework~\cite{BrezziF_CockburnB_MariniL_SuliE-2006aa}, the IP
methods for the more general case of mixed boundary conditions.  To
present a short derivation of the methods, we assume $\vect{u}\in
[H^{2}(\Omega)]^{d}$. Such assumption is not required for the methods
to work. We present the derivation under such assumption in order to
avoid unnecessary details which would shift the focus of our
presentation on preconditioners.

By assuming that the solution of~\eqref{eq:PDEform} is a priori
discontinuous, $\vect{u}\in [H^{2}(\Th)]^{d}$, we may rewrite the
continuous problem~\eqref{eq:PDEform} as follows: Find $\vect{u}\in
[H^{2}(\Th)]^{d}$ such that
\begin{equation}\label{wr00}
\left\{\begin{aligned}
-\div (\mathcal{C}\Gsym{u}) &=\vect{f} \quad &&\mbox{  on  } \K \in \Th\;, \\
\jump{ (\mathcal{C}\Gsym{u})\vect{n}}_E &=\vect{0} 
\quad &&\mbox{  on  } E \in \Eho\;,\\
\jump{\vect{u}}_E &=\vect{0} \quad &&\mbox{  on  } E \in \Eho\;,\\
\jump{\vect{u}}_E &=\vect{0} \quad &&\mbox{  on  } E \in \EhD\;,\\
\jump{ (\mathcal{C}\Gsym{u})\vect{n} -\vect{g}_N}_E &=\vect{0} 
\quad &&\mbox{  on  } E \in\EhN\;.
\end{aligned}\right.
\end{equation}
where we recall that $\mathcal{C}\Gsym{u} = 2\mu \Gsym{u} + \lambda
\operatorname{trace}(\Gsym{u}) I$.  Following
\cite{BrezziF_CockburnB_MariniL_SuliE-2006aa}, we next introduce a variational
formulation of \eqref{wr00} by considering the following five operators
\begin{align*}
\mathcal{B}_0: [H^{2}(\Th)]^{d}& \lor [L^{2}(\Th)]^{d},& 
\qquad\qquad \\
\mathcal{B}_1: [H^{2}(\Th)]^{d} &\lor [L^{2}(\Eho)]^{d},& \mathcal{B}^{\partial}_1: [H^{2}(\Th)]^{d} &\lor [L^{2}(\EhD)]^{d} \\
\mathcal{B}_2: [H^{2}(\Th)]^{d}& \lor [L^{2}(\Eho)]^{d},
& \mathcal{B}^{\partial}_2: [H^{2}(\Th)]^{d} &\lor [L^{2}(\EhN)]^{d},
\end{align*}
and weighting each equation in \eqref{wr00} appropriately. This then
amounts to considering the following problem: Find $\vect{u}\in
[H^{2}(\Th)]^{d}$ such that for all $\vect{v}\in [H^{2}(\Th)]^{d}$
\begin{equation}\label{wr0}
\begin{aligned}
  &(-\div(\mathcal{C}\Gsym{u})-\vect{f},
  \mathcal{B}_0(\vect{v}))_{\Th} + ( \jump{ (\mathcal{C}\Gsym{u})
    \vect{n}} ,
  \mathcal{B}_2(\vect{v}))_{\Eho} + ( \jump{\vect{u}}, \mathcal{B}_1(\vect{v}))_{\Eho}\qquad\qquad &&\\
  &\qquad\qquad \qquad\quad+ ( \jump{\vect{u}},
  \mathcal{B}^{\partial}_1(\vect{v}))_{\EhD}+ ( \jump{
    (\mathcal{C}\Gsym{u}) \vect{n}-\vect{g}_N} ,
  \mathcal{B}^{\partial}_2(\vect{v}))_{ \EhN}=\vect{0}.  &&
\end{aligned}
\end{equation}
Different choices of the operators $\calB_0$, $\calB_1$, $\calB_2$,
$\calB^{\partial}_1$ and $\calB^{\partial}_2$ above give rise to
different variational formulations and, consequently to different DG
methods. We refer to \cite[Theorem
6]{BrezziF_CockburnB_MariniL_SuliE-2006aa} for sufficient conditions
on the operators $\calB_0$, $\calB_1$, $\calB_2$, $\calB^{\partial}_1$
and $\calB^{\partial}_2$ to guarantee\footnote{We note that in
  \cite{BrezziF_CockburnB_MariniL_SuliE-2006aa} the focus is on the
  scalar Laplace equation. The arguments for the elasticity
  problem, are basically the same.}  the uniqueness of the solution of
\eqref{wr0}.

To derive the IP method of interest, we take $\vect{v}$ piecewise
smooth and we set $\mathcal{B}_0(\vect{v})=\vect{v}$,
$\mathcal{B}_2(\vect{v})=\avgv{v}$ and
$\mathcal{B}^{\partial}_2(\vect{v})=\vect{v}$ in \eqref{wr0}, to
obtain that
\begin{equation}\label{eq:ortho0}
\begin{array}{l}
(-\div (\mathcal{C}\Gsym{u}), \vect{v} )_{\Th} +
( \jump{( \mathcal{C}\Gsym{u})\vect{n}} , \avgv{v})_{\Eho\cup \EhD}+ 
(\jumpv{u}, \mathcal{B}_1(\vect{v}))_{\Eho\cup \EhD}\\
~~~= (\vect{f},\vect{v})_{\Th}+
( \vect{g}_N ,\vect{v} )_{\EhN}\;.
\end{array}
\end{equation}
Defining \begin{equation}\label{rhs00}
  \mathcal{F}(\vect{v})=(\vect{f},\vect{v})_{\Th}+( \jumpv{g} ,
  \mathcal{B}^{\partial}_1(\vect{v}) )_{\EhD} +( \vect{g}_N ,\vect{v}
  )_{\EhN},
\end{equation}
and integrating by parts the first term on the left side
of~\eqref{eq:ortho0} then leads to
\begin{equation}\label{eq:a-bilinear-form}
\scalarTT{\mathcal{C}\Gsym{u}}{\Gsym{v}} -
( \avg{(\mathcal{C}\Gsym{u})\vect{n}} , \jumpv{v} )_{\Eho\cup \EhD}+ 
(\jumpv{u},\mathcal{B}_1(\vect{v})) _{\Eho\cup \EhD}
= \mathcal{F}(\vect{v}).
\end{equation}
For a fixed edge $E\in\Eho\cup \EhD$ the operator $
\mathcal{B}_1(\vect{v})$ is defined by
\begin{equation}\label{defB1}
\mathcal{B}_1(\vect{v}):=
-\avg{(\mathcal{C}\Gsym{v}) \vect{n}} + 
\alpha_0\beta_0 \PzeroE\jumpv{v}+\alpha_1\beta_1\jumpv{v},
\end{equation}
where, following~\cite{HansboP_LarsonM-2003aa}, the parameters
$\beta_0$ and $\beta_1$ are chosen depending on the Lam\'{e} constants
$\lambda$ and $\mu$:
\begin{equation}\label{defbeta}
\beta_0:=d\lambda +2\mu, \qquad \quad \beta_1:=2\mu\;.
\end{equation}
The remaining two parameters, $\alpha_0$ and $\alpha_1$, are still at our
disposal to ensure (later on) stability and to avoid locking of the resulting method.
 
We define
\begin{equation}\label{eq:aj0aj1}
\begin{aligned}
a_{j,0}(\jumpv{u},\jumpv{v})&\de
\alpha_0\beta_0 \sum_{E\in \Eho\cup\EhD}
  \int_E\langle h_E^{-1}\jump{\vect{u}},
  \PzeroE\jumpv{v}\rangle,\\
 a_{j,1}(\jumpv{u},\jumpv{v})&\de 
\alpha_1\beta_1
  \sum_{E\in \Eho\cup\EhD}
  \int_E\langle h_E^{-1}\jump{\vect{u}},\jump{\vect{v}}\rangle\;,
\end{aligned}
\end{equation}
and set
\begin{equation*}
  a_j(\jumpv{u},\jumpv{v})=a_{j,0}(\jumpv{u},\jumpv{v})+  a_{j,1}(\jumpv{u},\jumpv{v}).
\end{equation*}
Then,  the weak formulation of Problem~\eqref{wr00} reads:
Find   $\vect{u}\in [H^{2}(\Th)]^{d}$ such that
\begin{equation}\label{eq:weakform-ok}
\calA(\vect{u},\vect{w})=\mathcal{F}(\vect{w}),\qquad \forall \,\,\vect{w}\in [H^{2}(\Th)]^{d}.
\end{equation}
The bilinear form $\calA(\cdot,\cdot)$ is given by
\begin{equation}\label{ipAstep}
\calA(\vect{u},\vect{w})=\calA_0(\vect{u},\vect{w})+ a_{j,1}(\jumpv{u},\jumpv{w}),
\end{equation}
where
\begin{equation}\label{ipA0}
\begin{array}{rcl}
\calA_0(\vect{u},\vect{w})& =& \scalarTT{\Ceps{u}}{\Gsym{w}}-
(\avg{(\Ceps{u})\vect{n}},\jumpv{w})_{\Eho\cup\EhD}\\
&& - ( \jumpv{u},\avg{(\Ceps{w})\vect{n}})_{\Eho\cup\EhD}
+ a_{j,0}(\jumpv{u},\jumpv{w}).
\end{array}
\end{equation}
It is straightforward to see that
\begin{equation}\label{ipA}
\begin{array}{rcl}
\calA(\vect{u},\vect{w}) &=& \scalarTT{\Ceps{u}}{\Gsym{w}}-
(\avg{(\Ceps{u})\vect{n}},\jumpv{w})_{\Eho\cup\EhD}\\
&& +\theta (\jumpv{u},\avg{(\Ceps{w})\vect{n}})_{\Eho\cup\EhD} + a_{j}(\jumpv{u},\jumpv{w}).
\end{array}
\end{equation}

To obtain the discrete formulation, we replace the function space 
$[H^{2}(\Th)]^{d}$ in ~\eqref{eq:weakform-ok} by $\Spacev{V}{DG}$, and
we get the {\bf IP-1 approximation} to the problem: Find $\uh\in \Spacev{V}{DG}$
such that:
\begin{equation}\label{eq:DG-formulation}
\calA(\uh,\vect{w})=\mathcal{F}(\vect{w}),\qquad \forall \,\,\vect{w}\in \Spacev{V}{DG}.
\end{equation}
We could also consider the approximation given by the {\bf IP-0
  method}: Find $\uh\in \Spacev{V}{DG}$ such that:
\begin{equation}\label{eq:DG-formulation0}
\calA_0(\uh,\vect{w})=\mathcal{F}(\vect{w}),\qquad \forall \,\,\vect{w}\in \Spacev{V}{DG}.
\end{equation}
As we see next, the {\bf IP-0} method provides a robust approximation
to the problem \eqref{eq:PDEform} in the pure displacement problem
$\Gamma_D=\partial\Omega$. As we mentioned earlier, for other types of
boundary conditions such equivalence in general does not hold.

\begin{remark}\label{theta}
  Although we do not consider non-symmetric IP methods in this paper,
  let us remark that non-symmetric versions can be easily incorporated
  in the definition of $\mathcal{B}_1(\vect{v})$.  For example, by
  setting:
\begin{equation*}
\mathcal{B}_1(\vect{v}):=
\theta\avg{(\mathcal{C}\Gsym{v}) \vect{n}} + 
\alpha_0\beta_0 \PzeroE\jumpv{v}+\alpha_1\beta_1\jumpv{v},
\end{equation*}
we obtain a non-symmetric bilinear form for the values $\theta=0$
or $\theta=1$. Such values of $\theta$ correspond to the Incomplete
Interior Penalty (IIPG, $\theta=0$) and Non-symmetric Interior Penalty
(NIPG, $\theta=1$) discretizations, respectively.
 
\end{remark}

\subsection{Stability Analysis\label{sec:2}} 
We close this section presenting the stability and continuity results pertinent to our work. We start by introducing some norm notation.
For $\vect{v}\in [H^2(\Th)]^d$ we define the semi-norms
\begin{equation}\label{semi-norm}
\begin{aligned}
&\|\nabla\vect{v}\|_{0,\Th}^{2}=\sum_{\K\in \Th}\|\nabla \vect{v}\|_{0,\K}^{2}  \qquad  && 
\|\mathcal{C}^{1/2}{\Gsym{v}}\|_{0,\Th}^{2}=\sum_{T\in\mathcal{T}_h}\int_{T}\langle\Ceps{v}:{\Gsym{v}}\rangle  \\ 
& | \calP^{0}_E\jump{\vect{v}}|^{2}_{\ast}=\sum_{E\in \Eho\cup\EhD}h_E^{-1} \|\calP^{0}_E\jumpv{v}\|_{0,E}^{2} \qquad &&|\jump{\vect{v}}|^{2}_{\ast}=\sum_{E\in \Eho\cup\EhD}h_E^{-1} \|\jumpv{v}\|_{0,E}^{2}\;,
\end{aligned}
\end{equation}
and norms:
\begin{equation}\label{norm_h}
\|\vect{v}\|_h^2 =
\|\mathcal{C}^{1/2}{\Gsym{v}}\|_{0,\Th}^{2} + \beta_0| \calP^{0}_E\jumpv{v}\|_{\ast}^{2}
+\beta_1|\jump{\vect{v}}|^{2}_{\ast}+ \sum_{E\in \Eho\cup\EhD} h_E\|{\bf \mathcal{C}}^{1/2} \Gsym{v}\cdot \n\|_{0,E}^{2}\;.
\end{equation} 
For  $\vect{v}\in \Spacev{V}{DG}$ we define the norms
\begin{equation}\label{eq:Me_D_norm}
\|\vect{v}\|_{DG0}^2=
\|\mathcal{C}^{1/2}{\Gsym{v}}\|_{0,\Th}^{2}
+\beta_0 | \calP^{0}_E\jump{\vect{v}}|^{2}_{\ast} 
\end{equation} 
and
\begin{equation}\label{normDG}
\|\vect{v}\|_{DG}^2 =\|\vect{v}\|_{DG0}^2+ \beta_1 
|\jump{\vect{v}}|^{2}_{\ast}\;.
\end{equation} 
Notice that for $\vect{v}\in \Spacev{V}{DG}$ the norms \eqref{norm_h} and \eqref{normDG} are equivalent. 
We finally introduce the norm:
\begin{equation}\label{fullnorm}
\|\vect{v}\|_{H^{1}(\Th)}^2 =\|\nabla \vect{v}\|_{0,\Th}^{2} + \beta_0| \calP^{0}_E\jumpv{v}\|_{\ast}^{2}++ \beta_1| \jumpv{v}\|_{\ast}^{2}\;.
\end{equation} 

Notice that continuity of the {\bf IP-1} and {\bf IP-0} bilinear forms
with respect to the norm \eqref{norm_h} follows easily from
Cauchy-Schwarz inequality together with the bound on the maximum
eigenvalue of $\mathcal{C}$, i.e., for all 
$\vect{u} \in [H^2(\Th)]^d$ and all $\vect{v} \in \Spacev{V}{DG}$ we have
\begin{eqnarray*}
  (\avg{(\Ceps{u}\vect{n})},\jump{\vect{v}})_{\Eho\cup\EhD}&=&  
  (\avg{(\Ceps{u}\vect{n})},\calP^{0}_E\jumpv{v})_{\Eho\cup\EhD}\\
  & \leq&   \frac{1}{\alpha_0\beta_0}  h^{1/2}_E\|{\bf \mathcal{C}} \Gsym{u}\cdot \n\|_{0,\Eho\cup\EhD} \cdot 
  \frac{  \alpha_0\beta_0}{4}  \|h_{E}^{-1/2}\calP^{0}_{E}\jump{\vect{v}}\|_{0,\Eho\cup\Gamma_D} \\
  & \leq &  \frac{1}{\alpha_0} \|h_{E}^{1/2}\mathcal{C}^{1/2}\eps{u}\cdot \n\|_{0,\Eho\cup\EhD}
  \frac{  \alpha_0\beta_0}{4}
  \|h_{E}^{-1/2}\calP^{0}_{E}\jump{\vect{v}}\|_{0,\Eho\cup\Gamma_D}.
\end{eqnarray*}
The equivalence of the norms \eqref{norm_h} and \eqref{normDG} for any
$\vect{v}\in \Spacev{V}{DG}$ guarantees therefore the continuity of the
{\bf IP-1} bilinear form with respect to the norm defined in
\eqref{normDG} for finite element functions.
 
The solvability of the discrete methods \eqref{eq:DG-formulation} and
\eqref{eq:DG-formulation0} is guaranteed if and only if, a discrete
version of the Korn's inequality holds on $\Spacev{V}{DG}$. In
\cite{BrennerS-2004aa} the following discrete Korn inequality is shown
for $[H^{1}(\Th)]^{d}$-vector fields:
\begin{equation}\label{korn:disc}
\|\vect{\nabla} \vect{v}\|_{0,\Th}^{2} \leq C \left( \|\Gsym{v}\|_{0,\Th}^{2} + |\pi_1 \jumpv{v}|_{\ast}^{2} +\|\vect{\nabla} \times \vect{v}\|_{0,\Th}^{2} \right) 
\end{equation}
where $\pi_1 : [L^{2}(\Eh)]^{d}\lor \vect{\mathbb{P}}^{1}(\Eh)$ is the
$L^{2}$-orthogonal projection onto the space of piecewise linear
vector valued functions on $\Eh$ (or a subset of it).

Coercivity of the {\bf IP-1} bilinear form with respect to the norm
\eqref{normDG} can be easily shown by taking
$\vect{u}=\vect{w}=\vect{v}$ in \eqref{ipA}:
\begin{align*}
\calA(\vect{v},\vect{v})
&=( \Ceps{\vect{v}} : \eps{v} )_{\Th}+ \alpha_0 \beta_0 \|h_{E}^{-1/2}\calP^{0}_{E}\jump{\vect{v}}\|_{0,\Eho\cup\Gamma_D}^{2} +
\alpha_1\beta_1\|h_{E}^{-1/2}\jump{\vect{v}}\|_{0,\Eho\cup\Gamma_D}^{2} &&\\
&\qquad \quad -2 (\avg{(\Ceps{v}\vect{n})}, \jump{\vect{v}})_{\Eho\cup\EhD}\;.
\end{align*}
Using Cauchy-Schwarz, trace and inverse inequalities together with the arithmetic-geometric inequality and the bound on the maximum eigenvalue
of $\mathcal{C}$ it follows that
\begin{eqnarray}\label{bella0}
(\avg{(\Ceps{v}\vect{n})},\jump{\vect{v}})_{\Eho\cup\EhD}&=&  
(\avg{(\Ceps{v}\vect{n})},\calP^{0}_E\jumpv{v})_{\Eho\cup\EhD} \nonumber \\
& \leq&   \frac{C_{t}(1+C_{inv})}{\alpha_0\beta_0} \|\mathcal{C}\eps{v}\|_{0,\Th}^{2} +
 \frac{  \alpha_0\beta_0}{4}  \|h_{E}^{-1/2}\calP^{0}_{E}\jump{\vect{v}}\|_{0,\Eho\cup\Gamma_D}^{2} \nonumber \\
  & \leq &  \frac{C_{t}(1+C_{inv})}{\alpha_0} \|\mathcal{C}^{1/2}\eps{v}\|_{0,\Th}^{2} +
\frac{  \alpha_0\beta_0}{4}  \|h_{E}^{-1/2}\calP^{0}_{E}\jump{\vect{v}}\|_{0,\Eho\cup\Gamma_D}^{2}.
\end{eqnarray}
Hence, we finally have
\begin{eqnarray*}
\calA(\vect{v},\vect{v})&\geq& 
(1-\frac{2C_{t}(1+C_{inv})}{\alpha_0} ) \|\mathcal{C}^{1/2}
\eps{v}\|_{0,\Th}^{2} +
\alpha_1\beta_1\|h_{E}^{-1/2}\jump{\vect{v}}\|_{0,\Eho\cup\Gamma_D}^{2}\\
&&~~~~+ \frac{\alpha_0}{2} \beta_0 \|h_{E}^{-1/2}\calP^{0}_{E}\jump{\vect{v}}\|_{0,\Eho\cup\Gamma_D}^{2}, \quad \forall\, \vect{v} \in \Spacev{V}{DG}\;,
\end{eqnarray*}
and therefore by taking $\alpha_{0}=\max{(1,4C_{t}(1+C_{inv}))}$
(sufficiently large) we ensure the coercivity of 
$\calA(\cdot,\cdot)$ with respect to the $\|\cdot\|_{DG}$-norm with constant
independent of $h$, $\mu$, and $\lambda$. Using now \eqref{korn:disc}
(since the norm \eqref{normDG} contains the full jump) we conclude
that $\calA(\cdot,\cdot)$ is coercive with respect to the
$\|\cdot\|_{H^{1}(\Th)}$-norm \eqref{fullnorm}. Therefore the {\bf
  IP-1} method defined by \eqref{ipA} provides a robust approximation
to \eqref{eq:PDEform} and does not lock as $\lambda \rightarrow
\infty$.

As we mentioned earlier, in the pure displacement case ($\Gamma_D=\partial\O,
\, \Gamma_N=\emptyset$) 
the bilinear form $\calA_0(\cdot\cdot)$ defined in
\eqref{ipA0} is coercive. Indeed we may use the identity (which holds
for $C_0^{\infty}(\Omega)$ functions):
\begin{equation}
\div \Gsym{v}=\frac{1}{2}\left( \div \nabla \vect{v} +\nabla \div \vect{v} \right)
\end{equation}
and rewrite the volume term in  \eqref{ipA0} (also in \eqref{ipA}) as
follows:
\begin{eqnarray*}
\scalarTT{\Ceps{u}}{\Gsym{w}}&=&
  \sum_{T\in\mathcal{T}_h}\int_{T}\langle\Ceps{v}:{\Gsym{v}}\rangle\\
&=& 
  \sum_{T\in\mathcal{T}_h}\int_{T} \left( 2\mu \langle \nabla
    \vect{u}:\nabla \vect{v}\rangle + 
(\mu +\lambda) \langle \div \vect{u}, \div \vect{v}\rangle \right).
\end{eqnarray*}
Then, from the discrete Poincar\'e inequality
\cite{FengX_KarakashianO-2001aa, BrennerS-2003aa}, the resulting
modified bilinear form for $\calA_0(\cdot,\cdot)$ is now coercive in
$\vect{V}^{DG}$ with respect to the $\|\cdot\|_{H^{1}(\Th)}$ norm,
with coercivity constant independent of $h$ and $\lambda$;
\begin{equation}\label{st:A00}
  \calA_0(\vect{v},\vect{v}) \geq C \|\vect{v}\|_{H^{1}(\Th)}^{2}  \quad \forall\, \vect{v}\in \vect{V}^{DG}\;.
\end{equation}
Therefore, the discrete problem \eqref{eq:DG-formulation0} is well
posed and the {\bf IP-0} method is stable and robust (locking free in
the limit $\lambda \to \infty$).  Notice that in \eqref{st:A00} we are
using the $\|\cdot \|_{H^{1}(\Th)}$-norm which includes not only the
norm $|P^{0}_E\jumpv{v}|_{\ast}$, but also the norm
$|\jumpv{v}|_{\ast}$. This is a consequence of the vector valued
counterpart of \cite[Lemma~2.3]{Ayuso-de-DiosB_ZikatanovL-2009aa}.
The stability property given in \eqref{st:A00} implies that the {\bf
  IP-0} and {\bf IP-1} methods are spectrally equivalent for the pure
displacement problem. These observations are summarized in the next
Lemma:
 \begin{lemma}\label{le:A0A1}
   Let $\mathcal{A}(\cdot,\cdot)$ and $\calA_0(\cdot,\cdot)$ be the
   bilinear forms of the {\bf IP-1} and {\bf IP-0} methods for the
   linear elasticity problem, defined in \eqref{ipA} and \eqref{ipA0},
   respectively. For the pure displacement problem
   $\Gamma_D=\partial\O, \,\, \Gamma_N=\emptyset$, there exist a
   constant $c>0$ that depends only on the geometry of the domain $\O$
   but is independent of the mesh size and the Lam\'e parameters $\mu$
   and $\lambda$ such that
 \begin{equation}\label{A0A1:0}
 \calA_0(\vect{v},\vect{v})\leq \calA (\vect{v},\vect{v})\leq c\calA_0(\vect{v},\vect{v}) \quad \forall\, \vect{v} \in \vect{V}^{DG}\;.
 \end{equation} 
 \end{lemma}
 The above lemma guarantees that for the pure displacement problem,
 constructing a uniform preconditioner for the {\bf IP-1} is
 equivalent to constructing a uniform preconditioner for the {\bf
   IP-0} method (see \cite{Ayuso-de-DiosB_ZikatanovL-2009aa}). For
 linear elasticity equations, unlike for scalar equations, this
 can be done only when $\Gamma_D=\partial\Omega$.

For a detailed derivation and error estimates, we refer
to~\cite[Theorem 2.5]{HansboP_LarsonM-2003aa}.

\section{Space decomposition}\label{sect:split}
We present now a decomposition of the DG space of piecewise linear
vector valued functions that plays a key role in the construction of
iterative solvers. This decomposition was introduced in
\cite{Ayuso-de-DiosB_ZikatanovL-2009aa} for scalar functions and also
in~\cite{BurmanE_StammB-2008aa} in a different context. Its extension
to vector valued functions is more or less straightforward.  We omit
those proofs which are just an easy modification of the corresponding
proofs in the scalar case. However, we review the main ingredients and
ideas behind such proofs, since they play an important role in the
analysis of the preconditioner given later on. In the last part of the
section we give some properties of the spaces entering in the and
prove a result that is essential for showing that the proposed
preconditioner is uniform.

Following \cite{Ayuso-de-DiosB_ZikatanovL-2009aa} we introduce the
space complementary to $\Space{V}{CR}$ in $\Space{V}{DG}$, 
 \begin{equation}\label{defZ}
\Space{\calZ}{}=\left\{ z\in \, \Space{V}{DG}~\mbox{and}~\calP_{E}^{0} 
\avg{z}=0,~\mbox{for all}~E\in \Eho \right\}.
\end{equation}
The corresponding space of vector valued functions is
\begin{equation}\label{defZ2}
\Spacev{\calZ}{}=[\Space{\calZ}{}]^d.
\end{equation}
To describe the basis functions associated with the spaces
\eqref{defCR} and \eqref{defZ2}, let $\varphi_{E ,\K}$ denote the
scalar basis function on $\K$, dual to the degree of freedom at the
mass center of the face $E$, and extended by zero outside $\K$. For
$E\in\partial\K$, $E^{\prime}\in \partial \K$, the function
$\varphi_{E,T}$ satisfies
\begin{equation*}
\varphi_{E,\K}(m_{E'})  = \left\{
\begin{array}{ll} 1 & \quad \mbox{if }\quad E=E^{\prime},\\ 
                0 &\quad \mbox{otherwise,} 
\end{array}
\right.
\end{equation*}
and also we have 
\begin{equation*}
\varphi_{E,\K}\in \mathbb{P}^{1}(\K), \quad
\varphi_{E,\K}(x)=0, \forall\, x\notin \K.
\end{equation*}
For all $\vect{u}\in \Spacev{V}{DG}$ we then have 
\begin{equation}\label{eq:represent0}
\vect{u}(x) = \sum_{\K\in \Th}\sum_{E\in \partial\K}\vect{u}_\K(m_E)\varphi_{E,\K}(x) =\sum_{E\in \Eh}
\vect{u}^{+}(m_E)\varphi_{E}^{+}(x)+\sum_{E\in \Eho}
\vect{u}^{-}(m_E)\varphi_{E}^{-}(x),
\end{equation}
where in the last identity we have just changed the order of
summation and used the short hand notation $\varphi_{E}^{\pm}(x):=\varphi_{E,\K^{\pm}}(x)$ together with
\begin{equation*}
\begin{aligned}
  \vect{u}^{\pm}(m_{E})&:=\vect{u}_{\K^{\pm}}(m_E)=\frac{1}{|E|}\int_{E}
\vect{u}_{\K^{\pm}} ds,
\quad &\forall\, E\in \Eho, \,\,:\,\, E=\partial\K^{+}\cap \partial \K^{-}, &&\\
  \vect{u}(m_{E})&:=\vect{u}_{\K}(m_E)=\frac{1}{|E|}\int_{E}
  \vect{u}_{\K} ds, \quad &\forall\, E\in \Ehb,~\mbox{such that}~E=\partial\K\cap \partial \O. &&
\end{aligned}
\end{equation*} 
Recalling now the definitions of $T^{+}(E)$ and $T^{-}(E)$ given
in~\eqref{eq:t-plus-minus} we set 
\begin{equation}\label{basis-CR}
\begin{aligned}
~& \varphi_{E}^{CR}=\varphi_{E,\K^{+}(E)}+\varphi_{E,\K^{-}(E)}, &\quad \forall\, E\in \Eho, &&\\
~& \varphi_{E}^{CR}=\varphi_{E,\K^{+}(E)}, &\quad \forall\, E\in \EhN. &&
\end{aligned}
\end{equation}
and
\begin{equation}\label{basis-Z}
\begin{aligned}
~&
\psi_{E}^{z}=\frac{\varphi_{E,\K^{+}(E)}-\varphi_{E,\K^{-}(E)}}{2},
&\quad \forall\, E\in \Eho, &&\\
~& \psi_{E}^{z}=\varphi_{E,\K^{+}(E)}, &\quad
\forall\, E\in \EhD.
\end{aligned}
\end{equation}
\begin{figure}[!hb]
\begin{center}
\includegraphics[height=2in,width=0.4\textwidth]{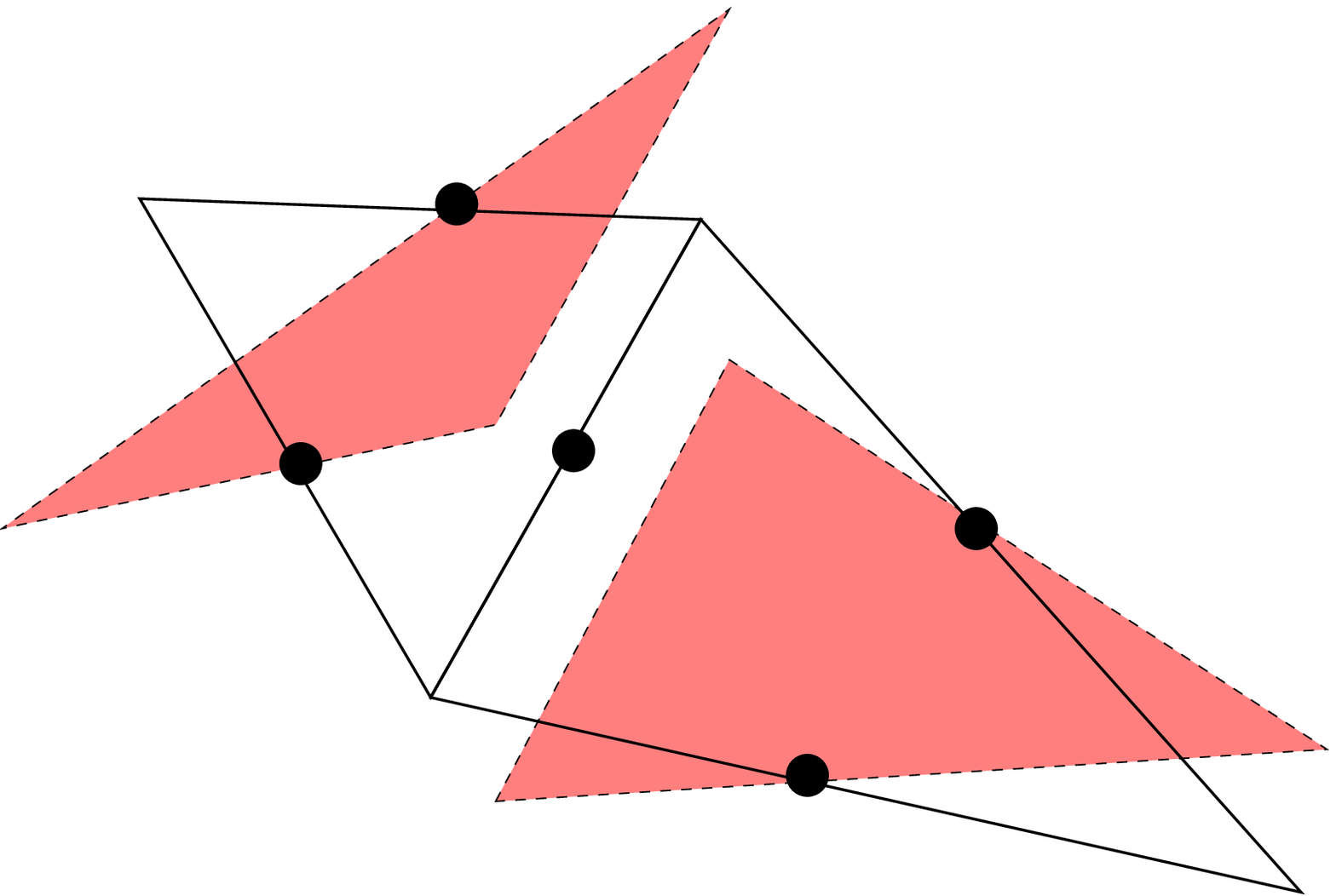}\hspace{0.025\textwidth}
\includegraphics[height=2in,width=0.4\textwidth]{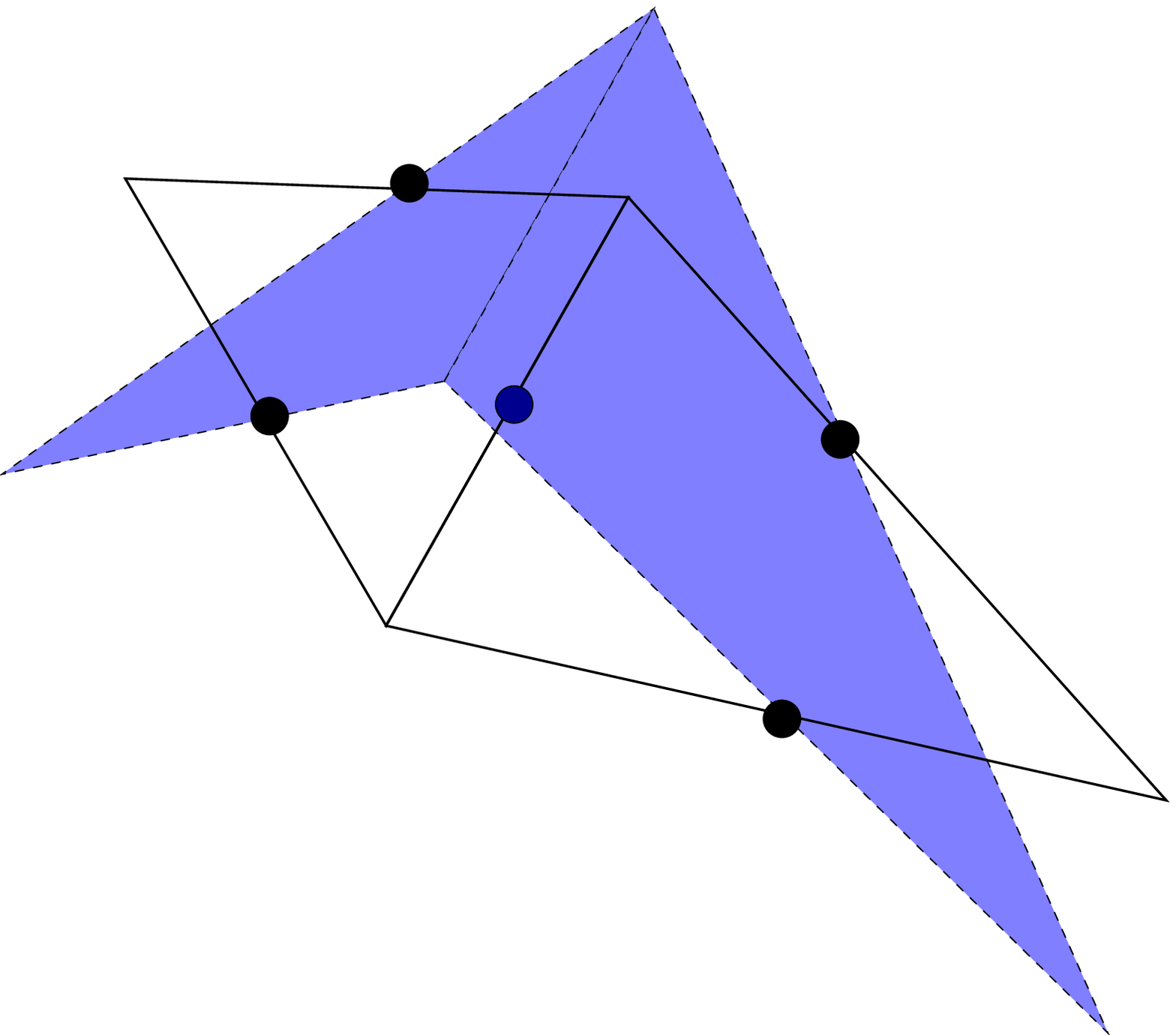}
\caption{Basis functions associated with the face $E$: $\psi_{E}^{z}$
  (left) and $\varphi_{E}^{CR}$ (right). \label{fig:ZCR}}
\end{center}
\end{figure}
Some clarification is needed here. Note that from the definition of
$\varphi_{E,\K^{+}(E)}$ and $\varphi_{E,\K^{-}(E)}$ for an interior
edge $E\in \Eho$, it does not follow that their sum is even defined on
the edge $E$, since it is just a sum of two functions from $L^2(\O)$.
However, the sum $(\varphi_{E,\K^{+}(E)}+\varphi_{E,\K^{-}(E)})$ has a
representative, which is continuous across $E$ and this representative
is denoted here with $\vect{ \varphi}_{E}^{CR}$, see Figure \ref{fig:ZCR}.

Clearly, $\{\varphi_{E}^{CR}\}_{E\in \Eho \cup\EhN}$ are
linearly independent, and
$\{\vect{\psi}^{z}_{E}\}_{E\in \Eho\cup\EhD}$
are linearly independent. A simple argument then
shows that 
\begin{equation*}
\Spacev{V}{CR} = \operatorname{span} \left\{
\{\varphi_{E}^{CR}\vect{e}_k\}_{k=1}^d\right\}_{E\in \Eho \cup\EhN}, \quad
\Spacev{\calZ}{} = \operatorname{span} 
\left\{
\{\psi_{E}^{z}\vect{e}_k\}_{k=1}^d\right\}_{E\in \Eho  \cup\EhD}.
\end{equation*}
Here $\vect{e}_k$, $k=1,\ldots,d$ is the $k$-th canonical basis vector in $\Reals{d}$. 
Hence by performing a change of basis in~\eqref{eq:represent0}, we
have obtained a ``natural'' splitting
of $$\Spacev{V}{DG}=\Spacev{V}{CR}\oplus \Spacev{\calZ}{}$$ and the set 
\begin{equation}\label{nat_basis}
 \{\vect{\psi}^{z}_{E}\}_{E\in \Eho\cup\EhD} \cup 
 \{\varphi_{E}^{CR}\}_{E\in \Eho \cup\EhN},
\end{equation} 
provides a natural basis for the DG finite element space.
This is summarized in
the next proposition.
\begin{proposition}\label{split:0}
For any $\vect{u}\in \Spacev{V}{DG}$
there exist unique $\vect{v}\in\Spacev{V}{CR}$ and a unique $\vect{z}\in \Spacev{\calZ}{}$ such that
\begin{equation}\label{eq:decomposition}
\begin{array}{ll}
\vect{u}=\vect{v}+\vect{z}\quad \mbox{and}\quad&
\begin{array}{rcl}
\vect{v}&=&\sum_{E\in   \Eho \cup\EhN}
\left(\frac{1}{|E|}\int_{E}\avgv{u}ds\right)\varphi^{CR}_{E}(x) \in \Spacev{V}{CR},\\
\vect{z}&=&\sum_{E\in   \Eho \cup\EhD}
\left(\frac{1}{|E|}\int_{E}\jump{\vect{u}}ds \right)\vect{\psi}^{z}_{E}(x) \in \Spacev{\calZ}{}.
\end{array}
\end{array}
\end{equation}
\end{proposition}
The proof of the above result follows by arguing as for the scalar
case in~\cite[Proposition 3.1]{Ayuso-de-DiosB_ZikatanovL-2009aa}, but
proceeding componentwise. The next Lemma shows that the splitting we have
proposed is orthogonal with respect to the inner product defined by
$\mathcal{A}_0(\cdot,\cdot)$.
\begin{lemma}\label{le:orto}
The splitting \eqref{eq:decomposition} $\Spacev{V}{DG}=\Spacev{V}{CR}\oplus  \Spacev{\calZ}{}$ is
$\calA_{0}$-orthogonal. That is
\begin{equation}\label{orto:1}
\calA_0(\vect{v},\vect{z})=\calA_0(\vect{z},\vect{v})=0 \qquad \forall\, \vect{v} \in \Spacev{V}{CR},\quad \forall\, \vect{z} \in \Spacev{\calZ}{} .
\end{equation}
\end{lemma}
The proof follows straightforwardly by using the weighted residual
formulation \eqref{eq:ortho0}-\eqref{ipA0} and the definition of the
spaces $ \Spacev{V}{CR}$ and $\Spacev{\calZ}{}$.

\subsection{Some properties of the space $\Spacev{\calZ}{}$\label{subsect:cbs}}
We now present some properties of the functions in the space
$\Spacev{\calZ}{}$. We start with a simple observation.  From the
definition of the spaces $\Spacev{V}{CR}$ and $\Spacev{\calZ}{}$ it is
easy to see that
\[
\sum_{\K\in \Th}\|\nabla \vect{z}\|_{0,\K}^2=(\jumpv{z},\avg{\nabla \vect{z}})_{\Eho\cup\EhD}.
\]
Applying the Schwarz inequality, one then gets the following estimate
\[
\sum_{\K\in \Th}\|\nabla \vect{z}\|^{2}_{0,\K}\leq
C\|h^{-1/2}P^{0}_E\jumpv{z}\|^{2}_{0,\Eh}, 
\]
which is a straightforward way to see that the restriction of the
{\bf IP-1} and {\bf IP-0}-bilinear forms (even for $\theta=0,1$ as in
Remark \ref{theta})  to the space $\Spacev{\calZ}{}$ are coercive in
the $\|\cdot\|_{H^{1}(\Th)}$-norm \eqref{fullnorm} (regardless whether
the boundary conditions are Dirichlet, Neumann or mixed type). Therefore the resulting stiffness matrices are positive definite.  

The next result provides bounds on the eigenvalues of
$\calA_0(\cdot,\cdot)$ and $\calA(\cdot,\cdot)$, when restricted to
$\Spacev{\calZ}{}$.

\begin{lemma}\label{le:poincZ}
Let $\Spacev{\calZ}{}$ be the space defined in \eqref{defZ2}. Then for all $\vect{z}\in
\Spacev{\calZ}{}$, the following estimates hold
\begin{equation}\label{poincZ:0}
  h^{-2}\|\vect{z}\|_{0}^{2}\lesssim \calA_{0}(\vect{z},\vect{z})\lesssim 
  h^{-2}\|\vect{z}\|_{0}^{2}\;,
\end{equation}
and also,
\begin{equation}\label{poincZ:1}
  [(\alpha_0)\beta_0 +\alpha_1\beta_1]  h^{-2}\|\vect{z}\|_{0}^{2}\lesssim \calA(\vect{z},\vect{z})\lesssim 
 [\alpha_0\beta_0 +\alpha_1\beta_1] h^{-2}\|\vect{z}\|_{0}^{2}\;,
\end{equation}
where $\beta_0$ and $\beta_1$ are as defined in  \eqref{defbeta}.
\end{lemma}
\begin{proof}
Arguing as in \cite[Lemma~5.3]{Ayuso-de-DiosB_ZikatanovL-2009aa} (but
now componentwise for vector valued functions) one can show  that
(due the special structure of the space $\Spacev{\calZ}{}$).
\begin{equation}\label{poinc:00}
  h^{-2}\|\vect{z}\|_{0}^{2}\lesssim \sum_{E\in \Eho\cup\EhD} h_{E}^{-1}\|\calP^{0}_{E}\jump{\vect{z}}\|^{2}_{0,E} \lesssim 
 h^{-2}\|\vect{z}\|_{0}^{2}\;.
 \end{equation}
From the coercivity of $\calA_0$  it follows then
 \begin{equation*}
   \alpha_0\beta_0 h^{-2}\|\vect{z}\|_{0}^{2}\lesssim \alpha_0\beta_0\sum_{E\in \Eho\cup\EhD} h_{E}^{-1}\|\calP^{0}_{E}\jump{\vect{z}}\|^{2}_{0,E} \leq \calA_0(\vect{z},\vect{z})\;.
   \end{equation*}
  Similarly, the $L^{2}(\Eh)$ stability of the projection $\calP_E^{0}$ together with the coercivity of $\calA$  gives
    \begin{equation*}
    \begin{aligned}
  ( \alpha_0\beta_0 +\alpha_1\beta_1) h^{-2}\|\vect{z}\|_{0}^{2}&\lesssim \alpha_0\beta_0\sum_{E\in \Eho\cup\EhD} h_{E}^{-1}\|\calP^{0}_{E}\jump{\vect{z}}\|^{2}_{0,E} + \alpha_1\beta_1\sum_{E\in \Eho\cup\EhD} h_{E}^{-1}\|\calP^{0}_{E}\jump{\vect{z}}\|^{2}_{0,E} &&\\
  &\lesssim \alpha_0\beta_0\sum_{E\in \Eho\cup\EhD}
  h_{E}^{-1}\|\calP^{0}_{E}\jump{\vect{z}}\|^{2}_{0,E} +
  C\alpha_1\beta_1\sum_{E\in \Eho\cup\EhD}
  h_{E}^{-1}\|\jump{\vect{z}}\|^{2}_{0,E} &&\\
& \leq 
   \calA(\vect{z},\vect{z}), &&
  \end{aligned}
   \end{equation*}
   and so, the lower bounds in \eqref{poincZ:0} and \eqref{poincZ:1}
   follow.  We next show the upper bound in \eqref{poincZ:0}, and the
   upper bound in \eqref{poincZ:1} is obtained in an analogous
   fashion. Using \eqref{bella0} together with \eqref{eigsC} we get
 \begin{align*}
\calA_0(\vect{z},\vect{z}) &\leq \alpha_0 \beta_0\sum_{E\in \Eho\cup \Gamma_D} h_{E}^{-1}\|\calP^{0}_{E}\jump{\vect{z}}\|^{2}_{0,E} + \|\mathcal{C}^{1/2}\eps{z}\|_{0,\Th}^{2} &&\\
&\leq \beta_0 \left(\alpha_0 \|h_{E}^{-1/2}\calP^{0}_{E}\jump{\vect{z}}\|^{2}_{\Eho\cup \Gamma_D} +C\|\eps{z}\|_{0,\Th}^{2}\right)\;.
\end{align*}
Hence, the upper bound in \eqref{poincZ:0} follows in a
straightforward fashion using the trace and inverse inequalities
together with the obvious inequality $\|\eps{z}\|_{0,\Th}\leq \|\grad
\vect{z}\|_{0,\Th}$.
\end{proof}

We close this section with establishing a uniform bound on the angle
between $\Spacev{V}{CR}$ and $\Spacev{Z}{}$ in the inner product given
by the bilinear form $\calA(\cdot,\cdot)$. The estimate is given in
Proposition~\ref{aux_aux_lema}. It plays a crucial role in bounding
the condition number of the preconditioned system.

We remind that $E \in\Eh$ denotes a $(d-1)$-dimensional simplex (a face),
which is either the intersection of two $d$-dimensional simplices $T\in \Th$
or an intersection of a $d$-dimensional simplex
$T\in\Th$ and the complement of $\Omega$, i.e., $E=T\cap
(\Reals{d}\setminus\Omega)$. In the former case, the face $E$ is called
an interior face and in the latter it is called a boundary face. 

The proof of Proposition~\ref{aux_aux_lema} requires arguments
involving the incidence relations between simplices $T\in\Th$ and
faces $E\in \Eh$, and estimates on the cardinality of these incidence
sets. For the readers' convenience, we provide a list of such
estimates below.

\begin{itemize}
\item We define $\mathcal{N}_0(E)$ to be the set of $d$-dimensional
$T\in \Th$  simplices that contain $E$:
\[
\mathcal{N}_0(E):=\{T\in \Th,\quad\mbox{such that}\quad E\in T\}
\]
By definition, for the cardinality of this set we have
$|\mathcal{N}_0(E)|=2$ for the interior faces and
$|\mathcal{N}_0(E)|=1$ for the boundary faces.

\item
We define the set of neighbor (or neighboring) faces $\mathcal{N}_1(E)$ to be
the set of faces which share an element with $E$:
\[
\mathcal{N}_1(E):=\{E'\in \Eh,\quad\mbox{such that}\quad 
\mathcal{N}_0(E)\cap\mathcal{N}_0(E')\neq \emptyset\}
\]
From Proposition~\ref{prop:cardinal} (see Appendix \ref{ap0}) we have that  $|\mathcal{N}_1(E)|\le (2d+1)$. 

\item Next, we define $\mathcal{N}_2(E)$ to be the set of faces
  which share at least one neighboring face with $E$:
\[
\mathcal{N}_2(E):=\{E'\in \Eh,\quad\mbox{such that}\quad
\mathcal{N}_1(E)\cap\mathcal{N}_1(E')\neq \emptyset\}
\]
From Proposition~\ref{prop:cardinal} we have the estimate
$|\mathcal{N}_2(E)|\le (2d+1)^2$.

\item For the basis functions~$\{\psi_E^z\}_{E\in \Eho\cup\EhD}$ we
  have the following relations: 
\begin{equation}\label{jumpZ}
\frac{1}{|E|}\int_E \jump{\psi_{E'}^z} = \delta_{EE'},\quad\mbox{and}\quad 
\jump{\psi_{E}^z}(x) = 1,\quad\mbox{for all}\quad x\in E, 
\end{equation}
\begin{equation}\label{jumpZ1}
|\jump{\psi_{E}^z}(x)| \le 1, \quad \mbox{for all}\quad x\in E',
\quad\mbox{and all}\quad E'\in \mathcal{N}_2(E).
\end{equation}
 The above relations all follow from the definition of
$\psi^z_E(x)$ and the fact that $\jump{\psi^z_E}$ is linear function on every face in $\Eh$,
and therefore 
$\int_E\jump{\psi^z_{E'}}=|E|\jump{\psi^z_{E'}}(m_{E})$. 

\item Finally, for 
$E\in \Eh$, $E'\in \Eh$, and $E''\in \Eh$ it is straightforward to see
that we have: 
\begin{equation}\label{jumpZ31}
\mbox{If}\quad E\notin\mathcal{N}_1(E')\cap  \mathcal{N}_1(E'')
\quad \mbox{then}\quad\int_E\jump{\psi^z_{E'}}\jump{\psi_{E''}^z}=0.
\end{equation}
An easy consequence from the definitions then is the following:
\begin{equation}\label{jumpZ3}
\mbox{If}\quad E'\notin  \mathcal{N}_2(E'')
\quad \mbox{then}\quad\int_E\jump{\psi^z_{E'}}\jump{\psi_{E''}^z}=0,
\quad \mbox{for all}\quad E\in \Eh.
\end{equation}
\end{itemize}
We finally give Proposition \ref{aux_aux_lema}. To avoid unnecessary 
complications with the notation, we state and prove the
result for scalar valued functions. The proof for vector valued functions is
easy to obtain, and with the same constant, 
by just applying the scalar valued result component-wise. 

\begin{proposition}\label{aux_aux_lema}
The following inequality holds for $z\in \Space{\calZ}{}$:
\begin{equation}\label{the_estimate}
\sum_{E\in \Eho\cup\EhD}\|h_E^{-1/2}(\jump{z}-\calP_{E}^{0}\jumpv{z})\|_{0,E}^2 
\le (1-\frac{1}{\rho})\sum_{E\in \Eho\cup\EhD}\|h_E^{-1/2}\jumpv{z}\|_{0,E}^2,
\end{equation}
with a constant $\rho\ge 1$ which depends on the shape regularity of the mesh.\end{proposition}
\begin{proof}

  Since $\calP_{E}^{0}$ is the $L^{2}$ orthogonal projection on the
  constants, we have that
\begin{equation}\label{ineq:zPz}
\|h_E^{-1/2}(\jump{z}-\calP_{E}^{0}\jumpv{z})\|_{0,E}^2 
=
\|h_E^{-1/2}\jump{z}\|_{0,E}^2- 
\|h_E^{-1/2}\calP_{E}^{0}\jumpv{z}\|_{0,E}^2. 
\end{equation}

Let $z \in \calZ$, i.e., $z=\sum_{E' \in \Eho\cup\EhD} z_{E'}
\psi_{E'}^{z}$. From~\eqref{jumpZ} we have that
$\calP_{E}^{0}\jump{\psi_{E'}^{z}} = \delta_{EE'}$, and hence, we may conclude that
\begin{eqnarray*}
\|h_E^{-1/2}\calP_{E}^{0}\jumpv{z}\|_{0,E}^2 &=& 
\sum_{E\in \Eho\cup\EhD}\sum_{E'\in\Eh}\delta_{EE'}\frac{|E|}{h_E}
z_Ez_{E'}\\
&=&\sum_{E\in \Eho\cup\EhD}\mathbb{D}_{EE} z_E^2=\langle \mathbb{D}\tilde{z},\tilde{z} \rangle. 
\end{eqnarray*}
Here we have denoted by $\mathbb{D}:\Reals{|\Eh|}\mapsto
\Reals{|\Eh|}$ a diagonal matrix with non-zero elements 
$\mathbb{D}_{EE}:=\frac{|E|}{h_E}$ and by 
$\tilde{z}\in \Reals{|\Eh|}$  the vector of coefficients
$\tilde{z}=\{z_E\}_{E\in  \Eh}$ in the expansion of $z\in \calZ$ via
the basis $\{\psi_{E}^{z}\}_{E\in \Eh}$.

Further we consider the right hand side of~\eqref{the_estimate} and we have 
\begin{eqnarray*}
\sum_{E\in \Eho\cup\EhD}\|h_E^{-1/2}\jumpv{z}\|_{0,E}^2 
&=& 
\sum_{E\in \Eho\cup\EhD}
h_E^{-1} \left\|\sum_{E' \in \Eho\cup\EhD} z_{E'} \jump{\psi_{E'}^{z}}\right\|_{0,E}^2\\
& = & \sum_{E\in \Eho\cup\EhD}
\int_E h_E^{-1} \sum_{E' \in \Eho\cup\EhD}\sum_{E'' \in \Eho\cup\EhD} z_{E'} z_{E''}
\jump{\psi_{E'}^{z}}\jump{\psi_{E''}^{z}}\\
&= &
\sum_{E' \in \Eho\cup\EhD}\sum_{E'' \in \Eho\cup\EhD} z_{E'} z_{E''}
\left(\sum_{E\in \Eh}\int_E h_E^{-1} 
\jump{\psi_{E'}^{z}}\jump{\psi_{E''}^{z}}\right)\\
&= &
\sum_{E' \in \Eho\cup\EhD}\sum_{E'' \in \Eho\cup\EhD} z_{E'} z_{E''}\mathbb{S}_{E'E''} = 
\langle \mathbb{S}\tilde{z},\tilde{z} \rangle. 
\end{eqnarray*}
Here, $\mathbb{S}:\Reals{|\Eh|}\mapsto \Reals{|\Eh|}$ denotes the
symmetric real matrix with elements
\begin{equation}\label{See}
\mathbb{S}_{E'E''}=
\sum_{E\in \Eho\cup\EhD} \int_Eh_E^{-1}\jump{\psi_{E'}^{z}}\jump{\psi_{E''}^{z}}=
\sum_{E\in\mathcal{N}_1(E')\cap \mathcal{N}_1(E'')}\int_Eh_E^{-1}\jump{\psi_{E'}^{z}}\jump{\psi_{E''}^{z}}.
\end{equation}
In the last identity above, we have used~\eqref{jumpZ31}. 
Note that
according to~\eqref{jumpZ3},  
if $E'\notin\mathcal{N}_2(E'')$ then $\mathbb{S}_{E'E''}=0$.
Thus,  
\[
\langle \mathbb{S}\tilde{z},\tilde{z} \rangle=
\sum_{E' \in \Eho\cup\EhD}\sum_{E''\in \mathcal{N}_2(E)} z_{E'} z_{E''}\mathbb{S}_{E'E''}.
\]
From this identity and~\eqref{jumpZ} and \eqref{jumpZ1}, we obtain
that 
\[
|\mathbb{S}_{E'E''}|\le 
|\mathcal{N}_1(E')\cap\mathcal{N}_1(E'')|\max_{E\in\mathcal{N}_1(E')\cap \mathcal{N}_1(E'')}
\frac{|E|}{h_E}\le
(2d+1)\max_{E\in\mathcal{N}_1(E')\cap \mathcal{N}_1(E'')}\frac{|E|}{h_E}.
\]
Introducing 
\[ 
\rho=\sup_{\tilde{w}\in \Reals{|\Eh|}}  \frac{\langle \mathbb{S} \tilde{w},\tilde{w}\rangle}{\langle \mathbb{D} \tilde{w},\tilde{w}\rangle},
=\sup_{\tilde{w}\in \Reals{|\Eh|}}  \frac{\langle \mathbb{D}^{-1/2}
  \mathbb{S} \mathbb{D}^{-1/2} 
\tilde{w},\tilde{w}\rangle}{\langle \tilde{w},\tilde{w}\rangle},
\] 
we obtain that
\begin{equation}\label{eq:ineq-s-d}
\langle \mathbb{S}\tilde{z},\tilde{z} \rangle = 
\langle
\mathbb{D}^{-1/2}\mathbb{S}\mathbb{D}^{-1/2}\mathbb{D}^{1/2}\tilde{z},
\mathbb{D}^{1/2}\tilde{z} \rangle \le
\rho \langle \mathbb{D}\tilde{z},\tilde{z} \rangle.
\end{equation}
This inequality can be rewritten as $\frac{1}{\rho}\langle \mathbb{S}\tilde{z},\tilde{z} \rangle \le
\langle \mathbb{D}\tilde{z},\tilde{z} \rangle$ and hence
\begin{equation*}
\langle \mathbb{S}\tilde{z},\tilde{z} \rangle -
\langle \mathbb{D}\tilde{z},\tilde{z} \rangle  \le 
\langle \mathbb{S}\tilde{z},\tilde{z} \rangle -
\frac{1}{\rho}\langle \mathbb{S}\tilde{z},\tilde{z} \rangle  = (1-\frac{1}{\rho})
\langle \mathbb{S}\tilde{z},\tilde{z} \rangle.
\end{equation*}

Note that~\eqref{ineq:zPz} implies that 
\begin{equation}\label{eq:rho-gt-1}
\langle \mathbb{S}\tilde{z},\tilde{z} \rangle =
\langle \mathbb{D}\tilde{z},\tilde{z} \rangle
+\sum_{E\in \Eh}\|h_E^{-1/2}(\jump{z}-\calP_{E}^{0}\jumpv{z})\|_{0,E}^2,
\end{equation}
and thus $\langle \mathbb{S}\tilde{z},\tilde{z} \rangle \ge \langle
\mathbb{D}\tilde{z},\tilde{z} \rangle$. This shows that $\rho
\ge 1$ in \eqref{eq:ineq-s-d}.  

It remains to show that $\rho$ can be bounded by quantities depending only on the shape
regularity of the mesh. Again, by~\eqref{jumpZ3} we have that:
if $E'\notin\mathcal{N}_2(E'')$ then $\mathbb{S}_{E'E''}=0$. Hence:
\begin{eqnarray*}
\rho &\le& \|\mathbb{D}^{-1/2}\mathbb{S}\mathbb{D}^{-1/2}\|_{\ell^{\infty}}
\le \max_{E'' \in \Eho\cup\EhD}\sum_{E'\in \mathcal{N}_2(E'')}
\frac{|\mathbb{S}_{E''E'}|}{\sqrt{\mathbb{D}_{E'E'}\mathbb{D}_{E''E''}}}\\
&\le& 
\max_{E'' \in \Eho\cup\EhD}\left[|\mathcal{N}_2(E'')|\max_{E'\in \mathcal{N}_2(E'')}
\frac{|\mathbb{S}_{E'E''}|}{\sqrt{\mathbb{D}_{E'E'}\mathbb{D}_{E''E''}}}\right]\\
&\le& (2d+1)^3\max_{E'' \in \Eho\cup\EhD}\max_{E'\in \mathcal{N}_2(E'')}
\max_{E\in\mathcal{N}_1(E')\cap \mathcal{N}_1(E'')}\frac{|E|}{h_E}
\sqrt{\frac{h_{E'}h_{E''}}{|E'||E''|}}.
\end{eqnarray*}
The quantity on the right side of this estimate only depends on the
 shape regularity of the mesh and the proof is complete. 
\end{proof}

\begin{remark}
  We remark that the constants in Proposition~\ref{aux_aux_lema} can
  be sharpened, at the price of further complicating the proof. The
  result given above is sufficient for our purposes, and we do not
  further comment on the possible ``optimal'' value of the constant
  $\rho$ above.  Another relevant observation is that the inequality
  in Proposition~\ref{aux_aux_lema} holds true, with the same or even
  smaller $\rho$, if we replace $\Eho\cup \EhD$ with a subset of edges
  $\mathcal{E}\subset (\Eho\cup \EhD)$ in~\eqref{the_estimate}. The
  proof is completely analogous (just $\Eho\cup \EhD$ is replaced by
  $\mathcal{E}$).
\end{remark}

\section{Preconditioning}\label{sect:preconditioning}

In this section, we present the construction and convergence analysis of the preconditioners we propose for the considered IP-methods. 

To construct the preconditioners, we use the subspace splitting given in Proposition \ref{split:0}, which suggests a simple change of basis. We have that for any $\vect{u},\vect{w}
\in \Spacev{V}{DG}$, we can write $\vect{u}=\vect{z}+\vect{v}$, and
$\vect{w}=\vect{\zeta}+\vect{\varphi}$, where $\vect{z},\vect{\zeta}
\in \Spacev{\calZ}{}$ and $\vect{v}, \vect{\varphi} \in
\Spacev{V}{CR}$.  Therefore, by performing this change of basis we can write 
$\mathcal{A}(\vect{u},\vect{w})=\mathcal{A}((\vect{z},\vect{v}),(\vect{\zeta},\vect{\phi}))$.
The $\calA_0$-orthogonality  \eqref{orto:1} of the subspaces in the splitting gives 
\[
\calA_0 ((\vect{z},\vect{v}),(\vect{\zeta},\vect{\phi})) = 
\calA_0 (\vect{z},\vect{\zeta}) +\calA_0 (\vect{v},\vect{\phi}).  
\]
which implies that the resulting stiffness matrix of $\calA_0$ in this new basis is block diagonal. 
For the pure displacement problem ($\Gamma_N=\emptyset$), as discussed
in Section \ref{sec:2}, the spectral equivalence given in Lemma
\ref{le:A0A1}, guarantees that an optimal preconditioner for $\calA_0$
is also  optimal for $\calA$. 
Therefore it is enough to study how to efficiently solve each of the blocks in the above block diagonal structure of $\calA_0$: the subproblem resulting from the restriction of  $\calA_0$ to $\Spacev{\calZ}{}$ and the subproblem on the space
$\Spacev{V}{CR}$.

For traction free or mixed type of boundary conditions, although a
preconditioner for $\calA_0$ does not result in an optimal solution
method. However, the block structure of $\calA_0$ in the new basis
already suggests that a reasonable choice for an approximation of
$\calA(\cdot,\cdot)$ is
\begin{equation}\label{ipB}
\calB ((\vect{z},\vect{v}),(\vect{\zeta},\vect{\phi})) = \calA (\vect{z},\vect{\zeta}) +\calA (\vect{v},\vect{\phi}).  
\end{equation}
The following algorithm describes the application of a
preconditioner, which is based on the bilinear form in the
equation~\eqref{ipB}. 
\begin{algorithm}\label{alg2}
  Let $\vect{r}\in [L^{2}(\O)]^d$ be given. Then the action of the
  preconditioner on $\vect{r}$ is the function $\vect{u}\in
  \Spacev{V}{DG}$ which is obtained from the following three steps.
  \renewcommand{\theenumi}{\arabic{enumi}.}
  \renewcommand{\labelenumi}{\theenumi~}
\begin{enumerate}
\item Find $z\in \Spacev{\calZ}{}$ such that  
\[
\calA(\vect{z},\vect{\zeta}) = (\vect{r},\vect{\zeta})_{\Th}
 \quad \mbox{for all}\quad \vect{\zeta}\in \Spacev{\calZ}{}.
\]
\item Find $\vect{v}\in \Spacev{V}{CR}$ such that 
\[
\calA(\vect{v},\vect{\varphi}) = (\vect{r},\vect{\varphi})_{\Th}
\quad \mbox{for all}\quad \vect{\varphi}\in \Spacev{V}{CR}.
\]
\item Set $\vect{u}=\vect{z}+\vect{v}$. 
\end{enumerate}
\end{algorithm}
As before, the application of this preconditioner corresponds
to solving 
the subproblem of the restriction of $\calA(\cdot,\cdot)$ to  $\Spacev{\calZ}{}$ and the subproblem of the restriction of  $\calA(\cdot,\cdot)$ to $\Spacev{V}{CR}$.

We now briefly discuss how the two smaller sub-problems can be
efficiently solved in both cases: (1) the case of Dirichlet boundary
conditions on all of $\partial\Omega$;  and (2) the
case of Neumann or mixed boundary conditions.

\noindent {\bf Solution in the subspace ${\bf \calZ}$:}
Lemma~\ref{le:poincZ} guarantees that the restriction of
$\calA(\cdot,\cdot)$ and $\calA_0(\cdot,\cdot)$ to $\Spacev{\calZ}{}$
is well-conditioned with respect to both, the mesh size and the
Lam\'{e} constants $\lambda, \mu$.  Therefore, the linear system
corresponding to the subproblem of the restriction to
$\Spacev{\calZ}{}$ can be efficiently solved by the method of
Conjugate Gradients (CG). A simple consequence of the well known
estimate on the convergence of CG (see, e.g., \cite{SaadY-2003aa,
  HestenesM_StiefelE-1952aa}) shows that the number of CG iterations
required to achieve a fixed error tolerance is uniformly bounded,
independently of the size of the problem and the parameters.

\noindent {\bf Solution in} $\Spacev{V}{CR}$: 

We now briefly discuss how to construct a uniform preconditioner for
the corresponding subproblem on the space $\Spacev{V}{CR}$.  Rather
than developing a completely new method, the idea is to use the
optimal preconditioners that have already been studied in literature,
and modify them if needed so that they fit in the present
framework. For our discussion, we distinguish two cases: the pure
displacement problem ($\Gamma_N=\emptyset$) and the case with mixed or
traction free boundary conditions ($\Gamma_N\ne \emptyset$).
  \begin{itemize}

  \item  For the case of Dirichlet boundary conditions on the entire
    boundary--the so-called pure displacement problem--it is known how
    to construct optimal order multilevel preconditioners that are
    robust with respect to the parameter $\lambda$, see
    e.g. \cite{BlahetaR_MargenovS_NeytchevaM-2006aa,KrausJ_MargenovS-2009aa,GeorgievI_KrausJ_MargenovS-2009aa}
    and the references therein.

  \item The traction free problem or the case of mixed boundary
    conditions is more difficult to handle because the (discrete) Korn
    inequality is not satisfied for the standard discretization by
    Crouzeix-Raviart elements without additional stabilization, as was
    shown in \cite{FalkR-1991aa}. The design of optimal and robust
    solution methods for stabilized discretizations is still an open
    problem, however, auxiliary space techniques might bridge this gap
    soon.
 \end{itemize}

\subsection{Convergence Analysis}
We now prove that the proposed block preconditioners are indeed optimal so that their convergence is uniform with respect to mesh size and the Lam\'e parameters. This result is given 
 in Theorem~\ref{the-theorem}.  The following Lemma is crucial for this
proof, since it gives estimates on the norm of the
off-diagonal blocks in the $2\times 2$ block form of the
stiffness matrix associated to $\mathcal{A}(\cdot,\cdot)$, corresponding to the space splitting
$\Spacev{V}{DG}=\Spacev{V}{CR}\oplus \Spacev{\calZ}{}$.  The result provides a measure of the angle between the subspaces $\Spacev{V}{CR}$ and $\Spacev{\calZ}{}$, with respect to the $\calA$-norm.
The proof of this result uses 
Proposition~\ref{aux_aux_lema}. 
\begin{lemma}\label{lem_gama}
{\bf Strengthened Cauchy-Schwarz inequality:}  
The following inequality holds for any $\vect{z}\in \Spacev{\calZ}{}$ and any
$\vect{v}\in \Spacev{V}{CR}$
$$
\calA(\vect{z},\vect{v})^2 \leq \gamma^2 \calA(\vect{z},\vect{z})
\calA(\vect{v},\vect{v})
$$
where $\gamma< 1$ and $\gamma$ depends only on $\alpha_0$, $\alpha_1$
and the constant from Proposition~\ref{aux_aux_lema}. 
\end{lemma}
\begin{proof}
  We know that we can always choose $\alpha_0$ large enough, such that
  for all $\vect{u} \in \Spacev{V}{DG}$ we have
\begin{equation*}
\calA_0(\vect{u},\vect{u})=\scalarTT{\Ceps{u}}{\Gsym{u}}
-2(\avg{(\Ceps{u}) \vect{n}},\jumpv{u})_{\Eho\cup\EhD}+ \alpha_0 a_{j,0}(\jumpv{u},\jumpv{u}) \ge 0.
\end{equation*}
Then it is sufficient to prove that there exists
$\gamma=\gamma(\alpha_1)<1$ such that for all
$\vect{z}\in \Spacev{\calZ}{}$ and for all $\vect{v} \in \Spacev{V}{CR}$ the inequality
$$
\left[a_{j,1}(\jumpv{z},\jumpv{v})\right]^2 \leq
\gamma^2 a_{j,1}(\jumpv{z},\jumpv{z}) a_{j,1}(\jumpv{v},\jumpv{v}),
$$
holds.  
By the definition of
the spaces $\Spacev{\calZ}{}$ and $\Spacev{V}{CR}$, on the boundary
edges $E\in \Ehb$ we have either 
$\calP_{E}^{0}\jump{\vect{z}}=0$ (if $E\in\EhN$) or
$\calP_{E}^{0}\jump{\vect{v}}=0$ (if $E\in\EhD$). Hence, from the symmetry of
$\calP_{E}^{0}$
we conclude that
\[
\int_E\langle \jump{\vect{z}},\calP_{E}^{0}\jump{\vect{v}}\rangle=
\int_E\langle
\calP_{E}^{0}\jump{\vect{z}},\jump{\vect{v}}\rangle=0, 
\quad \mbox{for all}\quad E\in \Ehb,
\quad \mbox{and all}\quad 
\vect{z}\in \Spacev{\calZ}{},
\quad\vect{v}\in 
\Spacev{V}{CR}.
\]
Since for the interior edges $E\in \Eho$ we also have
$\calP_{E}^{0}\jump{\vect{v}}=0$, the above relation and the
definition of $\PzeroE$ altogether imply that
for all
$\vect{z}\in \Spacev{\calZ}{}$,
and $\vect{v}\in 
\Spacev{V}{CR}$
\begin{equation}\label{eq:zero} 
\alpha_1\beta_1
  \sum_{E\in \Eho\cup\EhD}
  \int_E\langle h_E^{-1}\PzeroE\jump{\vect{z}},\jump{\vect{v}}\rangle=
 \alpha_1\beta_1
  \sum_{E\in \Eho\cup\EhD}
  \int_E\langle h_E^{-1}\jump{\vect{z}},\PzeroE\jump{\vect{v}}\rangle=0.
\end{equation} 
The equation~\eqref{eq:zero} and the Schwarz inequality then lead to
\begin{eqnarray*} 
[a_{j,1}(\jumpv{z},\jumpv{v})]^2 & = &
\left[  \alpha_1\beta_1
 \sum_{E\in \Eho\cup\EhD} \int_E\langle  h_E^{-1}\jumpv{z},\jump{\vect{v}}\rangle
 \right]^2 \\
&=&
\left[  \alpha_1\beta_1
 \sum_{E\in \Eho\cup\EhD} \int_E\langle
   h_E^{-1}(\jumpv{z}-\PzeroE\jump{\vect{z}}),\jump{\vect{v}}\rangle
 \right]^2 \\
& \leq &
a_{j,1}(\jumpv{v},\jumpv{v})  
\left[\alpha_1\beta_1\sum_{E\in \Eho\cup\EhD}
\|h_E^{-1/2}(\jumpv{z}-\PzeroE\jump{\vect{z}})\|_{0,E}^2\right].
\end{eqnarray*}
Next, the result in Proposition~\ref{aux_aux_lema} (more precisely its vector
valued form) implies that 
\begin{equation*}
\alpha_1\beta_1\sum_{E\in \Eho\cup\EhD}
\|h_E^{-1/2}(\jumpv{z}-\PzeroE\jump{\vect{z}})\|_{0,E}^2
 \le 
\left(1-\frac{1}{\rho}\right)  
\alpha_1\beta_1\sum_{E\in \Eho\cup\EhD} \|h_E^{-1/2}\jumpv{z}\|_{0,E}^2
\end{equation*}
Therefore, we have 
\begin{eqnarray*}
[a_{j,1}(\jumpv{z},\jumpv{v})]^2 
& \leq &
\left(1-\frac{1}{\rho}\right)  
  a_{j,1}(\jumpv{v},\jumpv{v})
\left[ \alpha_1\beta_1\sum_{E\in \Eho\cup\EhD} \|h_E^{-1/2}\jumpv{z}\|_{0,E}^2
 \right]\\
& \leq & \left(1-\frac{1}{\rho}\right)a_{j,1}(\jumpv{z},\jumpv{z}) a_{j,1}(\jumpv{v},\jumpv{v}),
\end{eqnarray*}
which shows the  desired inequality.
 \end{proof}

We are now in a position to prove that the preconditioner given by
Algorithm~\ref{alg2} is uniform with respect to the mesh size and the
problem parameters. 
\begin{theorem}\label{the-theorem}
  Let $\mathcal{A}(\cdot,\cdot)$ be the symmetric bilinear form
  defined by (\ref{ipA}) where $\theta=-1$ and $\calB(\cdot,\cdot)$ be
  the bilinear form defined by (\ref{ipB}). Then the following
  estimates hold for all $\vect{z}\in \Spacev{\calZ}{}$ and for all $\vect{v} \in \Spacev{V}{CR}$
\begin{equation}\label{eq:equiv-A-A_0}
\frac{1}{1+\gamma}\mathcal{A}((\vect{z},\vect{v}),(\vect{z},\vect{v}))
\leq
\calB((\vect{z},\vect{v}),(\vect{z},\vect{v}))
\leq
\frac{1}{1-\gamma}\mathcal{A}((\vect{z},\vect{v}),(\vect{z},\vect{v})).
\end{equation}
The constant $\gamma<1$ is the constant from Lemma~\ref{lem_gama}.
\end{theorem}
\begin{proof}
Using Lemma \ref{lem_gama} we have
$$
- 2 \gamma \sqrt{\calA(\vect{z},\vect{z}) \, \calA(\vect{v},\vect{v})}
\le
2 \calA(\vect{z},\vect{v})
\le
2 \gamma \sqrt{\calA(\vect{z},\vect{z}) \, \calA(\vect{v},\vect{v})}
$$
and since $-a^2-b^2 \le 2 a b \le a^2 + b^2$ for any real numbers $a$ and $b$ we obtain
$$
(1 - \gamma) \left( \calA(\vect{z},\vect{z}) + \calA(\vect{v},\vect{v}) \right)
\le
\calA(\vect{z},\vect{z}) + \calA(\vect{v},\vect{v}) + 2 \calA(\vect{z},\vect{v})
\le
(1 + \gamma) \left( \calA(\vect{z},\vect{z}) + \calA(\vect{v},\vect{v}) \right)
$$
which is the same as
$$
(1 - \gamma) \calB((\vect{z},\vect{v}),(\vect{z},\vect{v}))
\le
\calA((\vect{z},\vect{v}),(\vect{z},\vect{v}))
\le
(1 + \gamma) \calB((\vect{z},\vect{v}),(\vect{z},\vect{v}))
$$
and thus (\ref{eq:equiv-A-A_0}) holds with the same constant $\gamma < 1$ as used in the estimate of Lemma~\ref{lem_gama}.
\end{proof}

\begin{remark}\label{rem_gama} Note that  $\gamma\le q < 1$ is uniformly
  bounded away from $1$ and this bound holds independently of the
  parameters $h$, $\lambda$, and $\mu$.
\end{remark}

\section{Numerical experiments}\label{sect:numerics}

In this section we present a set of numerical tests that illustrate our theoretical results.
We consider the SIPG discretization of the model problem (\ref{eq:PDEform}) on the unit square in $\Reals{2}$
with mixed  boundary conditions. For the penalty parameters in (\ref{eq:aj0aj1}) we choose the values $\alpha_0=4$
and $\alpha_1=1$. The coarsest mesh (at level $0$) consists of eight triangles and is refined four times.
Each refined mesh at level $\ell$, $\ell=1,2,3,4$ is obtained by subdividing every triangle at level $(\ell-1)$ into
four congruent triangles. The CBS constants and the spectral condition numbers summarized in the tables below
have been computed using MATLAB.

In Table~\ref{gama_tabl} we list the values of the constant $\gamma^{2}$ in the inequality stated
in Lemma~\ref{lem_gama} for different levels of refinement. Evidently, $\gamma$ is uniformly bounded
with respect to the mesh size (or the number of refinement levels) and also with respect to the material parameters, 
Young's modulus $\mathfrak{E}$ and Poisson ratio $\nu$ (see Remark~\ref{rem_gama}).
\renewcommand{\arraystretch}{1}
\begin{table}[hb]
\begin{center}
 \caption{Observed CBS constant $\gamma^{2}$ for
    $\Omega=(0,1)^2$. \label{gama_tabl}}
\begin{tabular}{|c|ccccc|}
\hline
$\gamma^{2}$ & $\nu=0.25$ & $\nu=0.4$ & $\nu=0.49$ & $\nu=0.499$ & $\nu=0.49999$ \\
\hline
$\ell=1$ & 0.0664 & 0.025  & 0.0024 & 2.4024$\times10^{-4}$ & 2.4015$\times10^{-6}$ \\
$\ell=2$ & 0.0678 & 0.0255 & 0.0025 & 2.4567$\times10^{-4}$ & 2.4559$\times10^{-6}$ \\
$\ell=3$ & 0.0684 & 0.0258 & 0.0025 & 2.4866$\times10^{-4}$ & 2.4857$\times10^{-6}$ \\
$\ell=4$ & 0.0686 & 0.0259 & 0.0025 & 2.4974$\times10^{-4}$ & 2.4966$\times10^{-6}$ \\
\hline
\end{tabular}
\end{center}
\end{table}
It can be seen from Table~\ref{gama_tabl_jump} that the two subspaces $\Spacev{V}{CR}$ and
$\Spacev{\calZ}{}$ remain nearly $\calA$-orthogonal when we introduce a jump in the Poisson
ratio (on the coarsest mesh); In our experiment we set $\nu=\nu_1=0.3$ (and $E=E_1=1$) in
the subdomain $\Omega_1=[0,0.5]\times[0,0.5]\cup [0.5,1]\times[0.5,1]$, and $\nu=\nu_2$
(and $E_2=1$) in the subdomain $\Omega_2=\Omega \setminus \Omega_2$, respectively.
\begin{table}[hb]
\begin{center}
 \caption{Observed CBS constant $\gamma^{2}$ for $\Omega=(0,1)^2$ and
    jumps in $\nu$. \label{gama_tabl_jump}}
\begin{tabular}{|c|ccccc|}
\hline
$\gamma^{2}$ & $\nu_2=0.3$ & $\nu_2=0.4$ & $\nu_2=0.49$ & $\nu_2=0.499$ & $\nu_2=0.49999$ \\
\hline
$\ell=1$ & 0.0451 & 0.0177 & 0.0442 & 0.0509 & 0.0517 \\
$\ell=2$ & 0.0460 & 0.0180 & 0.0689 & 0.0803 & 0.0816 \\
$\ell=3$ & 0.0464 & 0.0182 & 0.0689 & 0.0802 & 0.0816 \\
$\ell=4$ & 0.0466 & 0.0182 & 0.0689 & 0.0802 & 0.0816 \\
\hline
\end{tabular}
\end{center}
\end{table}

Next we consider an L-shaped domain $\Omega=[0,1]\times[0,1] \setminus (0.5,1]\times(0.5,1]$
with Neumann boundary conditions on the sides $y=0$ and $y=1$ and Dirichlet boundary conditions
on the remaining part of the boundary. The initial triangulation (level 0) consists of 4 similar
triangles. The angle is almost the same as for the square domain, see Table~\ref{Lgama_tabl}.
\begin{table}[hb]
\begin{center}
\caption{Observed CBS constant $\gamma^{2}$ for L-shaped domain.\label{Lgama_tabl}}
\begin{tabular}{|c|ccccc|}
\hline
$\gamma^{2}$ & $\nu=0.25$ & $\nu=0.4$ & $\nu=0.49$ & $\nu=0.499$ & $\nu=0.49999$ \\
\hline
$\ell=1$ & 0.0561 & 0.0202 & 0.0019 & 1.8918$\times10^{-4}$ & 1.8906$\times10^{-6}$ \\
$\ell=2$ & 0.0631 & 0.0233 & 0.0022 & 2.2118$\times10^{-4}$ & 2.2106$\times10^{-6}$ \\
$\ell=3$ & 0.0672 & 0.0252 & 0.0024 & 2.4216$\times10^{-4}$ & 2.4207$\times10^{-6}$ \\
$\ell=4$ & 0.0682 & 0.0257 & 0.0025 & 2.4810$\times10^{-4}$ & 2.4801$\times10^{-6}$ \\
\hline
\end{tabular}
\end{center}
\end{table}

Furthermore, we computed the relative condition number of the preconditioner $B$ corresponding
to the bilinear form (\ref{ipB}) for the model problem on the L-shaped domain. The results of
this experiment, which are listed in Table~\ref{Lrel_cond_num}, confirm the uniform bound provided
by Theorem~\ref{the-theorem}.
\begin{table}
\begin{center}
\caption{Tabulated values of $\kappa(B^{-1}A)$ for L-shaped domain. 
\label{Lrel_cond_num}}
\begin{tabular}{|c|ccccc|}
\hline
$\kappa(B^{-1}A)$
 & $\nu=0.25$ & $\nu=0.4$ & $\nu=0.49$ & $\nu=0.499$ & $\nu=0.49999$ \\
\hline
$\ell=1$ & 1.6204 & 1.3314 & 1.0912 & 1.0279 & 1.0028 \\
$\ell=2$ & 1.6713 & 1.3606 & 1.0990 & 1.0302 & 1.0030 \\
$\ell=3$ & 1.6997 & 1.3774 & 1.1037 & 1.0316 & 1.0031 \\
$\ell=4$ & 1.7073 & 1.3820 & 1.1050 & 1.0320 & 1.0032 \\
\hline
\end{tabular}
\end{center}
\end{table}
   
Finally, we computed the condition number $\kappa(A_{zz})$ of the matrix $A_{zz}$ related to the restriction
of $\calA(\cdot,\cdot)$ to the space $\Spacev{\calZ}{}$, again for the model problem on the L-shaped domain.
In view of Lemma~\ref{le:poincZ} we already know that $A_{zz}$ is well-conditioned, and this is clearly seen in 
Table~\ref{Lcond_Z_block} where the values of $\kappa(A_{zz})$ are listed.
\begin{table}
\begin{center}
\caption{Values of $\kappa(A_{zz})$ for L-shaped
    domain. \label{Lcond_Z_block}}
\begin{tabular}{|c|ccccc|}
\hline
$\kappa(A_{zz})$ & $\nu=0.25$ & $\nu=0.4$ & $\nu=0.49$ & $\nu=0.499$ & $\nu=0.49999$ \\
\hline
$\ell=1$ & 8.9067 & 7.1484 & 6.4788 & 6.4220 & 6.4158 \\
$\ell=2$ & 9.0875 & 7.1932 & 6.4829 & 6.4229 & 6.4164 \\
$\ell=3$ & 9.1577 & 7.2080 & 6.4841 & 6.4230 & 6.4164 \\
$\ell=4$ & 9.1794 & 7.2118 & 6.4844 & 6.4230 & 6.4164 \\
\hline
\end{tabular}
\end{center}
\end{table}

\section{Acknowledgments}
Part of this work was completed while the fourth author was visiting
RICAM, Austrian Academy of Sciences in Linz. Thanks go to the RICAM
for the kind hospitality and support.  The work of the first author
was partially supported by the Spanish MEC under projects
MTM2008-03541 and HI2008-0173. The work of the second author has been
partially supported by the Bulgarian NSF, Grant DO 02-338/08. We also
gratefully acknowledge the support by the Austrian Science Fund,
Grants P19170-N18 and P22989-N18. The work of the fourth author has
been supported in part by the US National Science Foundation, Grants
DMS-0810982, and OCI-0749202.

\appendix
\section{Auxiliary results}\label{ap0}

\subsection{Bounds on the cardinality of $\mathcal{N}_1(E)$ and $\mathcal{N}_2(E)$} 
We first recall the definitions of $\mathcal{N}_0(E)$, $\mathcal{N}_1(E)$ and $\mathcal{N}_2(E)$, already
given in \S\ref{subsect:cbs}:
\begin{eqnarray*}
\mathcal{N}_0(E)&:=&\{T\in \Th,\quad\mbox{such that}\quad E\in T\},\\
\mathcal{N}_1(E)&:=&\{E'\in \Eh,\quad\mbox{such that}\quad 
\mathcal{N}_0(E)\cap\mathcal{N}_0(E')\neq \emptyset\},\\
\mathcal{N}_2(E)&:=&\{E'\in \Eh,\quad\mbox{such that}\quad 
\mathcal{N}_1(E)\cap\mathcal{N}_1(E')\neq \emptyset\}.
\end{eqnarray*}
In the proof of the strengthened Cauchy-Schwarz inequality
\S\ref{subsect:cbs} we needed several estimates on the cardinality of
these sets and these estimates are given in the proposition below. We
remind the reader that we have $|\mathcal{N}_0(E)|\le 2$.

\begin{proposition}\label{prop:cardinal}  
  The following inequalities hold:
\begin{equation}\label{eq:cardinal}
|\mathcal{N}_1(E)|\le (2d+1)\quad\mbox{and}\quad |\mathcal{N}_2(E)|\le (2d+1)^2.
\end{equation}
\end{proposition}
\begin{proof}
Let $E\in \Eh$ be fixed. To  prove the bound on
  $|\mathcal{N}_1(E)|$ we consider the elements
 $T\in \Th$, such that $E\in T$. In each such element $T$, there are exactly $d$
  faces $E'\in T$, $E'\neq E$. Since there are at most two elements $T\in
  \Th$ containing $E$ we have at most $2d$ faces $E'\in\Eh$ such
  that $E'\in \mathcal{N}_1(E)$, and $E'\neq E$. Adding $E$ itself to
  the total count gives $|\mathcal{N}_1(E)|\le (2d+1)$.

The second bound given in~\eqref{eq:cardinal} follows from the first
and the following inclusion:
\[
\displaystyle \mathcal{N}_2(E)\subset \bigcup_{E'\in \mathcal{N}_{1}(E)} \mathcal{N}_1(E').
\]
To show the above inclusion, we consider an arbitrary $E''\in
\mathcal{N}_2(E)$.  By the definition of $\mathcal{N}_2(E)$, the
intersection of $\mathcal{N}_1(E'')$ and $\mathcal{N}_1(E)$ is not
empty. Equivalently, there exists $E'\in\Eh$ such that $E'\in\mathcal
N_1(E'')$ and $E'\in\mathcal{N}_1(E)$.  On the other hand, from the
definition of $\mathcal{N}_1(E'')$, we have that
$E'\in\mathcal{N}_1(E'')$ implies that $E''\in\mathcal{N}_1(E')$,
i.e., if $E'$ is a neighbor of $E''$, then $E''$ is a neighbor of
$E'$.

Putting this together, we conclude that: if $E''\in \mathcal{N}_2(E)$,
then there exists $E'\in \mathcal{N}_1(E)$, such that
$E''\in\mathcal{N}_1(E')$, and this is exactly the inclusion we wanted
to show.

To prove the desired  bound is then straightforward:
\begin{eqnarray*}
\displaystyle \left|\bigcup_{E'\in \mathcal{N}_{1}(E)} \mathcal{N}_1(E')\right|&\le&
\sum_{E'\in \mathcal{N}_{1}(E)} |\mathcal{N}_1(E')|
 \le 
\sum_{E'\in \mathcal{N}_{1}(E)} (2d+1)
\\
&=&(2d+1)|\mathcal{N}_{1}(E)| 
\le 
(2d+1)^2.
\end{eqnarray*}
\end{proof}
\subsection{A multiplicative relation}\label{subsect:dumb}
This is to prove a basic relation used to derive~\eqref{eq:represent0}
as well as~\eqref{eq:ortho0}. 
Let $\odot$ be a map $V\times W\mapsto U$, where $U$, $V$, and $W$ are
linear vector spaces over the real numbers. We assume that $\odot$ satisfies the
following distributive laws:
\[
a\odot (b+c)= a\odot b + a\odot c, \qquad
(a+b)\odot c= a\odot c + b\odot c,
\]
and we assume that for all $\xi\in \Reals{}$
and all $\eta\in \Reals{}$, we have:
\begin{equation}\label{eq:prop0}
(\xi a)\odot (\eta b) = (\xi\eta) (a\odot b).
\end{equation}
We have the following identities, based on the
definitions~\eqref{eq:definition0}:
\begin{equation}\label{eq:identity0}
a^+\odot b^+-a^-\odot b^- = \jump{a}\odot \avg{b}+\avg{a}\odot \jump{b}.
\end{equation}
Proving this relation is indeed trivial. Some examples for which the
reader should verify these identities are: (1) For real numbers $a$
and $b$ one may take as $\odot$ the usual multiplication of real
numbers; (2) $a$ and $b$ elements of a real Hilbert space and $\odot$
inner product; (3) $a$ and $b$ are linear operators, and $\odot$ is
then the multiplication of linear operators. Note that in such case
$\odot$ is not necessarily commutative; (4) $a$ is a matrix and $b$ is
a vector, or more generally, $a$ is a linear operator and $b$ is an
element of a Hilbert space.

From~\eqref{eq:definition0}, we have that the right side of the
identity~\eqref{eq:identity0} is
\begin{equation*}
\jump{a}\odot \avg{b}+\avg{a}\odot \jump{b}=
(a^+ - a^-)\odot \left(\frac{b^+ + b^-}{2}\right)+
\left(\frac{a^+ + a^-}{2}\right)\odot (b^+ - b^-)
\end{equation*}
Using the distributive law, and~\eqref{eq:prop0} (linearity of $\odot$ with
respect to scalar multiplication), we have
\begin{eqnarray*}
\lefteqn{(a^+ - a^-)\odot \left(\frac{b^+ + b^-}{2}\right)+
\left(\frac{a^+ + a^-}{2}\right)\odot (b^+ - b^-)}\\
&=&
\frac12(a^+ - a^-)\odot (b^+ + b^-)+
\frac12(a^+ + a^-)\odot (b^+ - b^-)\\
&=&
\frac12a^+ \odot (b^+ + b^-)-
\frac12a^- \odot (b^+ + b^-)+
\frac12a^+ \odot (b^+ - b^-)+
\frac12a^- \odot(b^+ - b^-)\\
&=&
\frac12a^+ \odot b^+
+\frac12a^+ \odot b^-
-\frac12a^- \odot b^+
-\frac12a^- \odot b^-\\
&&
+\frac12a^+ \odot b^+ 
-\frac12a^+ \odot b^-
+\frac12a^- \odot b^+
-\frac12a^- \odot b^-\\
& = &
\frac12a^+ \odot b^+
-\frac12a^- \odot b^-
+\frac12a^+ \odot b^+ 
-\frac12a^- \odot b^-
=
a^+ \odot b^+ - a^- \odot b^-.
\end{eqnarray*}

\bibliographystyle{plain}

\end{document}